\documentclass[12pt]{amsart}

\usepackage{bbm} 
\usepackage{appendix}
\usepackage{enumerate}
\usepackage{tikz-cd}
\usepackage{color}
\usepackage[all]{xy}
\usepackage{amsmath, amssymb}
\usepackage{graphics}
\usepackage{graphicx} 
\usepackage{mathtools}
\usepackage[pdftex,hypertexnames=false,pagebackref=true]{hyperref}
\hypersetup{
	colorlinks=true,
	citecolor=blue,
	linkcolor=blue,
	urlcolor=blue
	}
	\usepackage{cleveref}
	\usepackage[square,numbers]{natbib}
	\usepackage{autonum}
\usepackage{geometry}
\numberwithin{equation}{subsection}
\theoremstyle{plain}
\newtheorem{theorem}{Theorem}[section]
\crefname{theorem}{Theorem}{Theorems}

\crefname{mainresult}{Main Result}{Main Results}
\newtheorem{lemma}[theorem]{Lemma}
\crefname{lemma}{Lemma}{Lemmas}
\newtheorem{proposition}[theorem]{Proposition}
\crefname{proposition}{Proposition}{Propositions}
\newtheorem{corollary}[theorem]{Corollary}
\crefname{corollary}{Corollary}{Corollaries}
\theoremstyle{definition}
\newtheorem{definition}[theorem]{Definition}
\crefname{definition}{Definition}{Definitions}
\newtheorem{example}[theorem]{Example}
\crefname{example}{Example}{Examples}
\newtheorem{remark}[theorem]{Remark}
\crefname{remark}{Remark}{Remarks}
\newtheorem*{UseOfAI}{Use of AI}

\newcommand{\Hom}{\operatorname{Hom}}

\renewcommand{\Im}{\operatorname{Im}}
\newcommand{\Ker}{\operatorname{Ker}}
\newcommand{\Coker}{\operatorname{Coker}\nolimits}
\newcommand{\id}{\operatorname{id}}

\newcommand{\rad}{\operatorname{rad}}

\newcommand{\Aut}{\operatorname{Aut}\nolimits}
\newcommand{\soc}{\operatorname{soc}\nolimits}

\newcommand{\supp}{\operatorname{supp}\nolimits}
\newcommand{\Irr}{\operatorname{Irr}}

\newcommand{\add}{\operatorname{add}}
\renewcommand{\mod}{\mathrm{mod}\,}

\newcommand{\Mod}{\mathrm{Mod}\,}

\newcommand{\proj}{\operatorname{proj}}

\newcommand{\Proj}{\mathcal{P}}

\newcommand{\Inj}{\mathcal{I}}
\newcommand{\Gproj}{\operatorname{Gproj}}

\newcommand{\sGproj}{\underline{\operatorname{Gproj}}\,}
\newcommand{\ind}{\operatorname{ind}}
\newcommand{\thick}{\operatorname{thick}}

\newcommand{\smod}{\underline{\operatorname{mod}}\,}

\newcommand{\F}{\operatorname{\mathbb{F}}\nolimits}
\newcommand{\Fp}{\mathbb{F}_\text{perf}}
\renewcommand{\P}{\operatorname{\mathbb{P}}\nolimits}
\newcommand{\Z}{\operatorname{\mathbb{Z}}\nolimits}

\newcommand{\gldim}{\operatorname{gl.dim}}
\newcommand{\injdim}{\operatorname{id}}

\renewcommand{\L}{\Lambda}
\newcommand{\G}{\Gamma}

\newcommand{\B}{\mathcal{B}}

\newcommand{\T}{\mathcal{T}}
\newcommand{\U}{\mathcal{U}}
\newcommand{\D}{\mathcal{D}}
\newcommand{\K}{\mathcal{K}}

\newcommand{\ZA}{\Z\hspace{-.6mm}A}
\newcommand{\Out}{\mathrm{Out}^{\rm tr}_K}

\newcommand{\A}{\mathcal{A}}
\newcommand{\C}{\mathcal{C}}
\newcommand{\sC}{\underline{\mathcal{C}}}

\newcommand{\X}{\mathcal{X}}
\newcommand{\op}{{\rm op}}

\newcommand{\Eco}{\mathbb{E}^{\rm co}}
\newcommand{\E}{\mathbb{E}}
\newcommand{\s}{\mathfrak{s}}
\newcommand{\tri}[6]{\!\!\xymatrix@C=5mm{#1\ar[r]^-{#2}&#3\ar[r]^-{#4}&#5\ar@{..>}[r]^-{#6}&}\!\!}

\newcommand{\csubstack}[2]{\left[\substack{\,#1\,\\[.5pt]\,#2\,}\right]}
\newcommand{\rsubstack}[2]{\left[\substack{\,#1\hspace{3pt}#2\,}\right]}

\newcommand{\xrightarrowtail}[1]{\stackrel{#1}{\rightarrowtail}}
\newcommand{\xtwoheadrightarrow}[1]{\stackrel{#1}{\twoheadrightarrow}}

\begin{document}
\title[G-radical-square-zero algebras]{G-radical-square-zero algebras}

\author[T. Honma]{Takahiro Honma}
\address{%
{\rm (T.H.)} Department of General Education \\
National Institute of Technology (KOSEN), Yuge College \\ 
Ochi \\
Ehime 794-2593\\ 
Japan}
\email{99cfqqc9@gmail.com}
\author[S. Usui]{Satoshi Usui}
\address{%
{\rm (S.U.)} Department of General Education \\
Tokyo Metropolitan College of Industrial Technology \\ 
8-17-1 Minamisenju, Arakawa-ku \\
Tokyo 116-8523\\ 
Japan}
\email{usui@metro-cit.ac.jp}

\thanks{%
T.H. was partially supported by JSPS KAKENHI Grant Number 26K16972.
S.U. was partially supported by JSPS KAKENHI Grant Number 25K17243.}

\keywords{%
G-radical-square-zero algebra,
Gorenstein-projective module,
Iwanaga-Gorenstein algebra,
singularity category,
triangulated category,
Jacobson radical,
Auslander-Reiten theory}

\subjclass[2020]{%
16G50, 
16G70, 
18A25
}

\date{}

\dedicatory{}

\begin{abstract}
In this paper, we first provide necessary and sufficient conditions for a Krull-Schmidt triangulated category with Auslander-Reiten triangles to have square-zero radical.
One of these conditions is that the middle term of every Auslander-Reiten triangle is either zero or indecomposable.
We then call an artin algebra {\it G-radical-square-zero} if the stable category of Gorenstein-projective modules has square-zero radical.
We characterize when a G-radical-square-zero algebra is CM-finite and, in this case, give an additive generator of the stable category.
For certain Iwanaga-Gorenstein G-radical-square-zero algebras, we further show that the stable category admits an additive generator consisting of higher syzygies of simple modules.
\end{abstract}

\maketitle
\tableofcontents


\section{Introduction}

This paper is motivated by the recent work of Hafezi and Bahlekeh \cite{Hafezi-Bahlekeh_2024}, who studied Gorenstein-projective  modules over certain artin algebras.
The notion of {\it Gorenstein-projective modules} was introduced by Enochs and Jenda \cite{Enochs-Jenda_1995_MathZ} and has its origins in the work of Auslander and Bridger \cite{Auslander-Bridger_1969} on {\it modules of Gorenstein dimension zero}.

The category $\Gproj\L$ of Gorenstein-projective modules over an artin algebra $\L$ is a Frobenius exact category whose projective objects are precisely the projective $\L$-modules.
Hence, its stable category $\sGproj\L$ carries the structure of a  triangulated category.
This stable category plays an important role in the representation theory of Iwanaga-Gorenstein algebras 
(i.e.~algebras $\L$ with  $\injdim_\L\L=\injdim_{\L^\op}\L<\infty$, where $\injdim_\L M$ denotes the injective dimension of a (right) $\L$-module $M$).
Indeed, a theorem of Buchweitz \cite{Buchweitz_1986} states that for an Iwanaga-Gorenstein algebra $\L$, $\sGproj\L$ is triangle equivalent to the {\it singularity category} $\D_{\rm sg}(\L)$ of $\L$, which is defined as the Verdier quotient of the bounded derived category $\D^{\rm b}(\mod \L)$ by the category $\K^{\rm b}(\proj \L)$ of perfect complexes. 
Note that an algebra $\L$ has finite global dimension if and only if its singularity category $\D_{\rm sg}(\L)$ vanishes.
Thus the singularity category measures the homological singularity of an algebra and hence captures homological information on an algebra of infinite global dimension. 
We also note that Gorenstein-projective modules are called {\it maximal Cohen-Macaulay modules} in Buchweitz's work.

Together with Buchweitz's theorem, the fact that Gorenstein-projective modules serve as projective objects in Gorenstein homological algebra highlights the importance of determining all indecomposable Gorenstein-projective modules over a given algebra. 
However, this task is usually difficult even in the case of CM-finite algebras (i.e.~algebras with only finitely many pairwise non-isomorphic indecomposable Gorenstein-projective modules).
Despite, many authors, such as \cite{Chen-Shen-Zhou_2018,Kalck_2015,Ringel_2013,Ringel-Zhang_2017}, obtained an explicit classification of indecomposable Gorenstein-projective modules for a specific class of algebras.
Among the algebras studied in these papers, all except those in \cite{Ringel-Zhang_2017} are CM-finite.

Recently, {\it G-semisimple algebras} were introduced by Hafezi and Bahlekeh \cite[Definition 3.1]{Hafezi-Bahlekeh_2024}.
As noted in \cref{remark_2}, an algebra $\L$ is G-semisimple if and only if the (Jacobson) radical $\rad\sGproj\L$ of $\sGproj\L$ vanishes  if and only if $\sGproj\L$ is a semisimple triangulated category.
It was proved in \cite[Proposition 3.4]{Hafezi-Bahlekeh_2024} that a G-semisimple algebra $\L$ is CM-finite if and only if $\Gproj\L$ is contravariantly finite in $\mod\L$.
This implies that Iwanaga-Gorenstein G-semisimple algebras are always CM-finite.
Moreover, it was proved in \cite[Proposition 3.13]{Hafezi-Bahlekeh_2024} that the following conditions are equivalent for an algebra $\L$:
\begin{enumerate}[(i)]
    \item $\L$ is an Iwanaga-Gorenstein G-semisimple algebra with $\injdim_\L\L\leq 1$.
    \item The radical $\rad\L$ of $\L$ additively generates $\sGproj\L$ (i.e.~$\sGproj\L=\add(\rad\L))$.
\end{enumerate}
In this case, all indecomposable non-projective Gorenstein-projective modules are obtained as direct summands of the first syzygy $\Omega_\L(\L/\rad\L)=\rad\L$ of the semisimple quotient $\L/\rad\L$.
This motivates us to find a higher dimensional analogue.

In this paper, we aim to extend the above results of Hafezi and Bahlekeh to the case of algebras $\L$ satisfying $\rad^2\sGproj\L=0$, which we call {\it G-radical-square-zero algebras}.
To this end, we first study triangulated categories $\T$ satisfying $\rad^2\T=0$, which we call {\it radical-square-zero triangulated categories}. 
Our first main result characterizes them in terms of Auslander-Reiten theory and the category $\mod\T$ of finitely presented contravariant 
functors from $\T$ to the category of abelian groups.

\begin{theorem}[{\cref{claim_3}}]
The following conditions are equivalent for a Krull-Schmidt triangulated category $\T$ with Auslander-Reiten triangles.
\begin{enumerate}
\item $\rad^2\T = 0$.
\item Every indecomposable projective object in $\mod\T$ is of length at most two.
\item The middle term of every Auslander-Reiten triangle in $\T$ is either zero or indecomposable.
\item Every connected component of the Auslander-Reiten quiver of $\T$ is of the form $\ZA_n/G$ for some $n \in\{1, 2\}$ and some weakly admissible group $G$ of automorphisms of $\ZA_n$.  
\end{enumerate}
\end{theorem}
Moreover, in \cref{claim_40,claim_36}, we discuss the triangulated structure of a radical-square-zero Krull-Schmidt triangulated category with Auslander-Reiten triangles.

We then study G-radical-square-zero algebras.
We generalize results for G-semisimple algebras in \cite[Section 3]{Hafezi-Bahlekeh_2024} to this setting.
Note that for any algebra $\L$, the condition $\rad^2\Gproj\L=0$ on the Frobenius category $\Gproj\L$ determines the structure of $\L$; see \cref{claim_22}.
We need some notations to state our second main result.
Denote by $X_M$ a minimal right $\Gproj\L$-approximation of a $\L$-module $M$, by $(-)^*=\Hom_\L(-, \L)$ the duality between $\sGproj\L$ and $\sGproj\L^\op$, by $\nu$ a Serre functor for $\sGproj\L$, and by $\thick(X)$ the smallest thick subcategory of $\sGproj\L$ containing an object $X$.

\begin{theorem}[{\cref{claim_19}}]\label{theorem_B}
The following conditions are equivalent for a G-radical-square-zero algebra $\L$.
\begin{enumerate}
    \item $\L$ is CM-finite.
    \item $\Gproj\L$ is contravariantly finite in $\mod\L$.
\end{enumerate}
In this case, 
\begin{align}
    \textstyle 
    \add\left(X_{\rad\L}\oplus \Omega_{\L^\op}(X_{\rad\L_{\L^\op}})^*\right) 
    = \sGproj\L 
    = \add\left(X_{\rad\L} \oplus \nu X_{\rad\L}\right).
\end{align}
In particular, 
\begin{align}
    \textstyle 
    \sGproj\L  =     
    \begin{cases}
    \add\left(X_{\rad\L}\right) & \text{if $\L$ is G-semisimple,}\\[5pt]
    \add\left(X_{\rad\L}\oplus \Omega_\L^{|l|}(X_{\rad\L})\right) & \text{if $\sGproj\L $ is $l$-Calabi-Yau.}
    \end{cases}
\end{align}
In both cases, $\sGproj\L = \thick(X_{\rad\L})$.
\end{theorem} 
The theorem recovers the results \cite[Proposition 3.4]{Hafezi-Bahlekeh_2024} and \cite[Proposition 3.13(i)$\Rightarrow$(ii)]{Hafezi-Bahlekeh_2024} mentioned above.
We extend \cite[Proposition 3.13]{Hafezi-Bahlekeh_2024} itself in \cref{claim_35}.  
In case $\L$ is a finite dimensional CM-finite G-radical-square-zero algebra over an algebraically closed field, one can determine whether $\sGproj\L$ is Calabi-Yau as explained in \cref{remark_3} (cf.~\cref{example_8}).

Specializing \cref{{theorem_B}} to Iwanaga-Gorenstein algebras, we obtain \cref{claim_20}, where we show that if $\L$ is an Iwanaga-Gorenstein G-radical-square-zero algebra with $\injdim_\L\L\leq d$, then it is CM-finite, and we may replace $X_{\rad\L}$ and $X_{\rad\L_{\L^\op}}$ with the $d$-th syzygies $\Omega_\L^d(\L/\rad\L)$ and $\Omega_{\L^\op}^d(\L/\rad\L)$ of $\L/\rad\L$, respectively.
Therefore, $\sGproj\L$ has a nice additive generator in the following cases:
\begin{align}
    \sGproj\L = 
    \begin{cases}
    \add\left(\Omega_\L^d(\L/\rad\L)\right) &\text{if $\L$ is G-semisimple,}\\[1.5mm]
    \add\left(\Omega_\L^d(\L/\rad\L)\oplus\Omega_\L^{d+|l|}(\L/\rad\L)\right) &\text{if $\sGproj\L$ is $l$-Calabi-Yau.}
    \end{cases}
\end{align}
In the former case, computing the $d$-th syzygies of all simple $\L$-modules leads to determining all indecomposable non-projective Gorenstein-projective $\L$-modules.
A similar conclusion can be obtained in the latter case as well.
We give examples of these two cases within the class of monomial algebras; see \cref{example_4,example_3}.

Applying Buchweitz's equivalence \cite{Buchweitz_1986}, we obtain the following description of the singularity category $\D_{\rm sg}(\L)$ of an Iwanaga-Gorenstein algebra $\L$ satisfying $\rad^2\D_{\rm sg}(\L)=0$.

\begin{theorem}[{\cref{claim_41}}]\label{theorem_C}
Let $\L$ be an Iwanaga-Gorenstein  G-radical-square-zero algebra.
Then 
\begin{align}
    \textstyle 
    \add\left(\L/\rad\L \oplus \Sigma^{1+n_{\L}+n_{\L^\op}}(\L/\rad\L)^*\right)=\D_{\rm sg}(\L)
    = \add\left(\L/\rad\L \oplus \nu\L/\rad\L\right),
\end{align}
where $\Sigma$ denotes the shift functor of $\D_{\rm sg}(\L)$, and $n_{\L}$ and $n_{\L^\op}$ are defined in \eqref{eq_13}.
In particular, 
\begin{align}
    \textstyle 
    \D_{\rm sg}(\L)=     
    \begin{cases}
    \add\left(\L/\rad\L\right) & \text{if $\D_{\rm sg}(\L)$ is semisimple,} \\[1pt]
    \add\left(\L/\rad\L \oplus \Sigma^l\L/\rad\L\right) & \text{if $\D_{\rm sg}(\L)$ is $l$-Calabi-Yau.}
    \end{cases}
\end{align}
In both cases,  $\D_{\rm sg}(\L) = \thick(\L/\rad\L)$.
\end{theorem}

This paper is organized as follows.
In \Cref{Preliminaries}, we recall some definitions and facts used in this paper.
In \Cref{Triangulated_categories_with_radical_square_zero}, we characterize the vanishing of the radical square of a Krull-Schmidt triangulated category with Auslander-Reiten triangles.
In \Cref{G-radical_square_zero_algebras}, we apply results in \Cref{Triangulated_categories_with_radical_square_zero} to Gorenstein-projective modules and extends results in \cite[Section 3]{Hafezi-Bahlekeh_2024} to the case of G-radical-square-zero algebras.
In Appendix \ref{Semisimple_triangulated_categories}, we provide necessary and sufficient conditions for a Krull-Schmidt triangulated category to be semisimple.

\begin{UseOfAI}
GPT-6 was used to improve the English writing and readability of this paper.
\end{UseOfAI}



\section{Preliminaries}\label{Preliminaries}

We first fix some conventions and terminologies.
Throughout this paper, we denote by $R$ (resp.~$\mathbbm{k}, k, K$) a commutative ring (resp.~a commutative artin ring, a field, an algebraically closed field).
By a module, we mean a finitely generated right module.
For a ring $\L$, we denote by $\mod\L$ the category of $\L$-modules and by $\proj\L$ the category of projective $\L$-modules.
We denote by 
$\injdim_{\L}M$ the injective dimension of a $\L$-module $M$.
We denote by $|X|$ the cardinality of a set $X$.
We assume that all categories are essentially small, and all subcategories are full and closed under isomorphisms.
Let $\C$ be an additive category.
We write the morphism set $\Hom_\C(X, Y)$ as $\C(X, Y)$.
For a collection $\X$ of objects in $\C$, we denote by $\add\X$ the smallest subcategory of $\C$ which is closed under finite direct sums and direct summands and contains $\X$. 
If $\X=\{X\}$ is a singleton, we write $\add(X) := \add\X$.
We say that $\C$ is {\it $R$-linear} if all morphism sets are $R$-modules and the composition of morphisms is $R$-bilinear.
We say that an $R$-linear category $\C$ is {\it Hom-finite} if the morphism set $\C(X, Y)$ has finite length over $R$ for all $X, Y \in \C$.  
We assume that all functors between $R$-linear categories are $R$-linear (i.e.~all induced morphisms between morphism sets are morphisms of $R$-modules). 
Finally, we denote by $\Sigma$ the shift functor of a triangulated category.

\subsection{Sink and source morphisms}
\label{Sink_and_source_morphisms}

In this subsection, we recall the definitions of sink and source morphisms and some related facts that will be used in this paper.

Let $f: X \to Y$ be a morphism in an additive category $\C$.
We say that $f: X \to Y$ is {\it right minimal} if every $h: X \to X$ such that $fh =f$ is an isomorphism.
We say that $f: X \to Y$ is {\it right almost split} if it is a non-retraction, and any non-retraction $g: W \to Y$ factors through $f: X \to Y$.
A non-zero object $X \in \C$ is called {\it indecomposable} if $X_1 =0$ or $X_2 =0$ whenever there is an isomorphism $X \cong X_1 \oplus X_2$.
Observe that if $f: X \to Y$ is right almost split, then $Y$ is indecomposable. 
A {\it sink morphism} is a morphism that is both right minimal and right almost split. 
Dually, we define a {\it left minimal morphism}, a {\it left almost split morphism} and a {\it source morphism}.
If there exists a sink morphism $f:X \to Y$, then it is unique in the sense that if there exists another sink morphism $f^\prime:X^\prime \to Y$, then $f=f^\prime g$ for some isomorphism $g:X \to X^\prime $. Dually for a source morphism.
For an indecomposable object $Y\in\C$, we denote by $u_Y: S_Y \to Y$ (resp.~$u^Y: Y\to S^Y$) a sink morphism ending (resp.~a source morphism starting) at $Y$ if it exists.
A non-zero object $X \in \C$ is called {\it endo-local} if its endomorphism ring $\C(X, X)$ is local. 
Observe that any endo-local object is indecomposable (\cite[Lemma 4.5]{Shah_2023}).
We say that {\it $\C$ has sink morphisms} if every endo-local object $Y\in\C$ has a sink morphism ending at $Y$. 
Dually, we define the condition that {\it $\C$ has source morphisms}.
It is clear that $\C$ has sink morphisms if and only if $\C^\op$ has source morphisms.

A non-zero object $X$ in an abelian category $\A$ is called {\it simple} if its only subobjects are $0$ and $X$.
A non-zero object $X \in \A$ is called {\it semisimple} if it is isomorphic to a (possibly infinite) coproduct of simple objects.

We denote by $\Mod\C$ the category of contravariant additive functors from an additive category $\C$ to the category $Ab$ of abelian groups.
A morphism in $\Mod\C$ is a natural transformation.
It is known that $\Mod\C$ is an abelian category (\cite[Chapter 5]{Freyd_book}).
We say that $F\in\Mod\C$ is {\it finitely generated} if there exists an exact sequence
\begin{align}
  \C(-,A) \to F \to 0
\end{align}
in $\Mod\C$ for some $A\in \C$.
We say that $F\in\Mod\C$ is {\it finitely presented} if there exists an exact sequence
\begin{align}
    \C(-,B) \to \C(-,A) \to F \to 0
\end{align}
in $\Mod\C$ for some $A, B\in \C$.
We denote by $\mod\C$ the subcategory of $\Mod\C$ consisting of finitely presented functors.
Recall that a {\it weak kernel} of a morphism $f: Y \to Z$ is a morphism $g: X \to Y$ such that the induced sequence
\begin{align}
    \C(-,X) \xrightarrow{(-,g)} \C(-,Y) \xrightarrow{(-,f)} \C(-,Z) 
\end{align}
in $\Mod\C$ is exact.
We say that {\it $\C$ has weak kernels} if every morphism admits a weak kernel. 
Dually, we define a {\it weak cokernel} of $f$ and the condition that {\it $\C$ has weak cokernels}.
It is known that $\mod\C$ is an abelian category if and only if $\C$ has weak kernels (\cite[Theorem 1.4]{Freyd_1966}).

In the rest of this section, let $\C$ be a Krull-Schmidt category \cite{Krause_2015}, that is, an additive category such that every object $X \in \C$ has a finite direct sum decomposition $X=X_1  \oplus \cdots \oplus X_n$ into endo-local objects.
It is known that every morphism in a Krull-Schmidt category is of the form $[0\ m]: A \oplus B \to C$, where $m: B \to C$ is right minimal (\cite[Lemma 1]{Rump_2005}).
We denote by $\ind\C$ the set of isomorphism classes of indecomposable objects in $\C$.
Recall that the {\it (Jacobson) radical} $\rad\C$ of $\C$ is the ideal  given by 
\begin{align}
    \rad\C(X, Y) =\{ h \in \C(X, Y) \mid 1_X - gh \text{ is invertible for any } g \in \C(Y, X) \}
\end{align}
for $X, Y\in\C$.
Let $\C/\rad\C$ be the ideal quotient of $\C$ by $\rad\C$.
This category has the same objects as $\C$, and its morphisms are given by 
\begin{align}
    (\C/\rad\C)(X, Y):= \C(X, Y)/\rad\C(X, Y)
\end{align}
for any $X, Y \in \C/\rad\C$.
By \cite[Proposition 2.3]{Auslander_1974_II} (cf.~\cite[Lemma A.5]{Krause_2012}), there is a bijection between $\ind\C$ and the set of isomorphism classes of simple functors in $\Mod\C$, which assigns to each $Z\in\ind\C$ the functor $(\C/\rad\C)(-, Z)$.
Moreover, by \cite[Corollary 2.6]{Auslander_1974_II}, an indecomposable object $Z\in\ind\C$ admits a right almost split morphism $v: Y \to Z$ if and only if $(\C/\rad\C)(-, Z) \in \mod\C$. 
In this case, the right almost split morphism $v: Y \to Z$ induces the following projective presentation of $(\C/\rad\C)(-, Z)$:
\begin{align}
    \C(-,Y) \xrightarrow{(-, v)} \C(-,Z) \xrightarrow{\pi} (\C/\rad\C)(-, Z) \to 0,
\end{align}
where $\pi$ is the canonical epimorphism, which is a projective cover of $(\C/\rad\C)(-, Z)$ (\cite[Lemma A.3]{Krause_2012}).
Moreover, $v: Y \to Z$ is a sink morphism if and only if the above presentation is minimal (\cite[Lemma A.2]{Krause_2012}).

A morphism $f\in\C(X, Y)$ is called  {\it irreducible} if it is neither a section nor a retraction, and either $g$ is a retraction or $h$ is a section for any factorization $f=gh$.
If $X, Y \in\ind\C$, then $f\in\C(X, Y)$ is irreducible if and only if $f\in \rad\C(X, Y)\setminus\rad^2\C(X, Y)$, where $\rad^2\C$ stands for the ideal of $\C$ consisting of morphisms of the form $gh$ with $g, h\in\rad\C$. 
For any $X, Y \in \ind\C$, 
we define
\begin{align}
    (\rad\C/\rad^2\C)(X, Y):= \rad\C(X, Y)/\rad^{2}\C(X, Y)
\end{align}
and set 
\begin{align}\label{eq_2}
d_{X,Y}:=\dim\,(\rad\C/\rad^2\C)(X, Y)_{(\C/\rad\C)(X,X)}
\end{align}
and
\begin{align}\label{eq_3}
d_{X,Y}^\prime := \dim  {}_{(\C/\rad\C)(Y,Y)}(\rad\C/\rad^2\C)(X, Y).
\end{align}
Although the standard notation for $(\rad\C/\rad^2\C)(X, Y)$ is $\Irr_\C(X,Y)$, we do not use it in this paper.
The following fundamental result plays a crucial role in our investigation.

\begin{proposition}[{\cite[Proposition 3.12]{Iyama-Nakaoka-Palu_2024}}]
Let $\C$ be a Krull-Schmidt category.
\begin{enumerate}
\item Every sink morphism $u_Z:S_Z \to  Z$ in $\C$ induces an isomorphism 
\begin{align}
    (-, u_Z) : (\C/\rad\C)(-, S_Z) \xrightarrow{\sim} (\rad\C/\rad^2\C)(-,Z)
\end{align}
of functors. 
In particular, we have $S_Z = \bigoplus_{W\in\ind\C}W^{d_{W,Z}}$ and thus $d_{W,X}<\infty$ for all $W\in\ind\C$.
\item Every source morphism $u^X:X \to S^X$ in $\C$ induces an isomorphism 
\begin{align}
    (u^X, -) : (\C/\rad\C)(S^X, -) \xrightarrow{\sim} (\rad\C/\rad^2\C)(X, -)
\end{align}
of functors. 
In particular, we have $S^X = \bigoplus_{W\in\ind\C}W^{d_{X,W}^\prime}$ and thus $d_{X,W}^\prime<\infty$ for all $W\in\ind\C$.
\end{enumerate}    
\end{proposition}

From the proposition, if $Z\in\ind\C$ admits a sink morphism $u_Z:S_Z\to Z$, then for any indecomposable summand $X$ of $S_Z$ with canonical inclusion $\iota: X \hookrightarrow S_Z$, the composition $u_Z \circ \iota: X \to Z$ is an irreducible morphism.
Dually for a source morphism.



\subsection{Extriangulated categories} 
\label{Extriangulated categories}

In this subsection, we briefly recall the definition of an extriangulated category and describe Auslander-Reiten theory for extriangulated categories.
We refer to \cite{Iyama-Nakaoka-Palu_2024,Nakaoka-Palu_2019_extriangulated} for more details.

Two sequences $A \xrightarrow{x} B\xrightarrow{y} C$ and $A \xrightarrow{x^\prime} B^\prime\xrightarrow{y^\prime} C$ in an additive category $\C$ are called {\it equivalent} if these fit into a commutative diagram
\begin{align}
    \xymatrix@M=1mm{
    A \ar[r]^-{x} \ar@{=}[d] & B \ar[r]^-{y} \ar[d]^-{b}& C \ar@{=}[d]\\
    A \ar[r]^-{x^\prime} & B^\prime \ar[r]^-{y^\prime} &C\\ 
    }
\end{align}
for some isomorphism $b:B\to B^\prime$.
The equivalence class of $A \xrightarrow{x} B\xrightarrow{y} C$ is denoted by $[A \xrightarrow{x} B\xrightarrow{y} C]$.
We say that $\C$ is {\it extriangulated} if it is equipped with a functor $\E:\C^\op \times \C \to Ab$ and a correspondence $\s$ assigning to each $\delta\in\E(C,A)$ an equivalence class $[A \to B \to C]$ of sequences in $\C$, where $\E$ and $\s$ satisfy the axioms $(\mathrm{ET}1)$, $(\mathrm{ET}2)$, $(\mathrm{ET}3)$, $(\mathrm{ET}3)^\op$,  $(\mathrm{ET}4)$ and $(\mathrm{ET}4)^\op$ given in \cite[Definition 2.12]{Nakaoka-Palu_2019_extriangulated}.

In \cite[Example 2.13 and Proposition 3.22]{Nakaoka-Palu_2019_extriangulated}, Nakaoka and Palu showed that exact and triangulated categories are fundamental examples of extriangulated categories.
They also characterized exact and triangulated categories among extriangulated categories; see \cite[Corollaries 3.18 and 7.6]{Nakaoka-Palu_2019_extriangulated}.
Note that Msapato \cite[Theorem 3.3]{Msapato_2025} characterized triangulated extriangulated categories in terms of inflations and deflations, which are defined below.
In what follows, we regard exact and triangulated categories as extriangulated categories in the sense of Nakaoka and Palu.

In the rest of this subsection, let $(\C, \E, \s)$ satisfy $(\mathrm{ET}1)$, $(\mathrm{ET}2)$, $(\mathrm{ET}3)$ and $(\mathrm{ET}3)^\op$.
A sequence $A \xrightarrow{x} B \xrightarrow{y} C$ in $\C$ is called a {\it conflation} if  $\s(\delta)=[A \xrightarrow{x} B \xrightarrow{y} C]$ for some $\delta \in \E(C, A)$.
A morphism $x : A\to B$ is called an {\it inflation} if it appears in a conflation $A \xrightarrow{x} B \xrightarrow{} C$.
A morphism $y : B\to C$ is called a {\it deflation} if it appears in a conflation  $A \xrightarrow{} B \xrightarrow{y} C$.
A pair of a conflation $A \xrightarrow{x} B \xrightarrow{y} C$ and an element $\delta \in \E(C, A)$ satisfying $\s(\delta)=[A \xrightarrow{x} B \xrightarrow{y} C]$ is called an {\it $\E$-triangle}.
We denote such a pair by $\tri{A}{x}{B}{y}{C}{\delta}$.
We say that an $\E$-triangle $\tri{X}{u}{Y}{v}{Z}{\delta}$ {\it splits} if $\delta=0$.
We say that $X\in\C$ is {\it projective} (resp.~{\it injective}) if  $\E(X, \C)=0$ (resp.~$\E(\C, X)=0$). 
We denote by $\Proj(\C) = \Proj_\E(\C)$ (resp.~$\Inj(\C)=\Inj_\E(\C)$) the subcategory of projective (resp.~{\it injective}) objects.
For a triangulated category $\T$, we have $\Proj(\T)=\{0\}=\Inj(\T)$, since $\T(-, \Sigma-)$ plays a role of $\E$.

For any $\delta\in\E(C,A)$ and any $a \in\C(A,A^\prime)$ and $c\in\C(C^\prime, C)$, we define 
\begin{align}
a\delta:=\E(C, a)(\delta)\in\E(C, A^\prime)
\quad\text{and}\quad
\delta c :=\E(c,A)(\delta)\in\E(C^\prime, A),
\end{align}
where $\E(C, a):  \E(C, A) \to  \E(C, A^\prime)$ and $\E(c, A):  \E(C, A) \to  \E(C^\prime, A)$.
Note that $(a\delta) c = \E(c, a)(\delta) = a(\delta c)$, since $\E(c, A^\prime) \E(C, a)  = \E(c, a) = \E(C^\prime, a) \E(c, A)$.

Following \cite[Definition 2.1 and Proposition 2.10]{Iyama-Nakaoka-Palu_2024},  
we say that a non-split $\E$-triangle $\tri{X}{u}{Y}{v}{Z}{\delta}$ is an {\it Auslander-Reiten $\E$-triangle} if $u$ is left almost split and $v$ is right almost split. 
In this case, we have $u=u^X$ and $v=u^Z$; see \cite[Theorem 2.9]{Iyama-Nakaoka-Palu_2024}.
By \cite[Proposition 2.6]{Iyama-Nakaoka-Palu_2024}, for any two Auslander-Reiten $\E$-triangles $\tri{X}{}{Y}{}{Z}{\delta}$ and $\tri{X^\prime}{}{Y^\prime}{}{Z^\prime}{\delta^\prime}$, the following conditions are equivalent: 
\begin{enumerate}
    \item There are isomorphisms $x\in\C(X, X^\prime)$ and $z\in\C(Z^\prime, Z)$ such that $x \delta z = \delta^\prime$.
    \item $X\cong X^\prime$.
    \item $Z\cong Z^\prime$.
\end{enumerate}
For this reason, we write $\tau Z:= X$ and $\tau^{-1} X:= Z$.

We say that {\it $\C$ has right Auslander-Reiten $\E$-triangles} if for any endo-local non-projective object $Z\in\C$, there exists an Auslander-Reiten $\E$-triangle $\tri{X}{}{Y}{}{Z}{}$ ending at $Z$. 
Dually, we define the condition that {\it $\C$ has left Auslander-Reiten $\E$-triangles}.
We say that {\it $\C$ has Auslander-Reiten $\E$-triangles} if it has right and left Auslander-Reiten $\E$-triangles.
By \cite[Lemma 3.2]{Iyama-Nakaoka-Palu_2024}, if $\C$ is a Krull-Schmidt category with sink (resp.~source) morphisms, then it has right (resp.~left) Auslander-Reiten $\E$-triangles.
The converse holds true when $\C$ is triangulated (since $\Proj(\C)=\{0\}=\Inj(\C)$).

We call a conflation $\delta:X \xrightarrowtail{u} Y\xtwoheadrightarrow{v}Z$ in an exact category an {\it Auslander-Reiten conflation} if the associated $\E$-triangle $\tri{X}{u}{Y}{v}{Z}{\delta}$ is Auslander-Reiten. 
An exact category is said to {\it have right (resp.~left) Auslander-Reiten conflations} if it has right (resp.~left) Auslander-Reiten $\E$-triangles as an extriangulated category.
An exact category is said to {\it have  Auslander-Reiten conflations} if it has right and left Auslander-Reiten conflations.

A triangle $X\xrightarrow{u}Y\xrightarrow{v}Z\xrightarrow{w}\Sigma X$ in a triangulated category is Auslander-Reiten in the usual sense of Happel \cite[4.1]{Happel_book} if and only if the associated $\E$-triangle $\tri{X}{u}{Y}{v}{Z}{w}$ is Auslander-Reiten.
A triangulated category is said to {\it have right (resp.~left) Auslander-Reiten triangles} if it has right (resp.~left) Auslander-Reiten $\E$-triangles as an extriangulated category.
A triangulated category is said to {\it have Auslander-Reiten triangles} if it has right and left Auslander-Reiten triangles.





\section{Triangulated categories with radical square zero} 
\label{Triangulated_categories_with_radical_square_zero}

The aim of this section is to characterize radical-square-zero triangulated categories via Auslander–Reiten theory.
Moreover, we investigate the triangulated structure of such a category $\T$, using the notion of path-closed subsets introduced by Chen and Ringel \cite{Chen-Ringel_2018}. 
Throughout this section, let
$\C$ be a Krull-Schmidt category and $\T$ a Krull-Schmidt triangulated category.

Define the {\it support} of a functor $F \in \Mod\C$ as $\supp F := \{ X \in \ind\C \mid F(X) \neq 0 \}$.
Following \cite[Theorem 2.12]{Auslander_1974_II}, we say that $\Mod\C$ is {\it locally finite} if for each object $X \in\C$, $|\supp\C(-, X)| < \infty$, and the right $\C(Y, Y)$-module $\C(Y, X)$ has finite length for each $Y \in \ind\C$.
See \cite[Proposition 1.11]{Auslander_1974_II} for necessary and sufficient conditions for $\Mod\C$ to be locally finite.

\begin{remark}
The dual of Propositions \ref{claim_16}, \ref{claim_29} to \ref{claim_39} below also hold. 
We leave the details to the reader.
\end{remark}

\begin{proposition}[Auslander] \label{claim_16}
If $\Mod\C$ is locally finite, then  $\C$ has sink morphisms. 
\end{proposition}
\begin{proof}
If $\Mod\C$ is locally finite, then all simple functors are finitely presented by \cite[Proposition 1.11 and Theorem 2.12]{Auslander_1974_II}.
Thus \cite[Proposition 2.7]{Auslander_1974_II} concludes that $\C$ has sink morphisms.
\end{proof}

\begin{example}\label{example_1}
Assume that $\C$ is $R$-linear Hom-finite with $|\ind\C|<\infty$.
Then $\Mod\C$ and $\Mod(\C^\op)$ are both locally finite, and hence  $\C$ has both sink morphisms and source morphisms by \cref{claim_16}.
\end{example}

\begin{example}\label{example_2}
We say that $\T$ is {\it locally finite} if $\Mod\T$ and $\Mod(\T^\op)$ are locally finite.
It follows from \cref{claim_16} and \Cref{Extriangulated categories} that if $\T$ is locally finite, then it has both sink morphisms and source morphisms, or equivalently, it has Auslander-Reiten triangles.
This was already obtained in \cite[Corollary 3.2]{Krause_2012}.
\end{example}

\begin{definition}
An additive category $\D$ is called {\it radical-square-zero} if $\rad^2\D=0$.  
\end{definition}

The converse in \cref{claim_16} holds true for a radical-square-zero Krull-Schmidt category.

\begin{proposition}\label{claim_29}
Assume that $\rad^2\C=0$.
Then $\Mod\C$ is locally finite if and only if $\C$ has sink morphisms.
\end{proposition}
\begin{proof}
By \cref{claim_16}, it suffices to show the \lq\lq if\,\rq\rq\,part.
Assume that $\C$ has sink morphisms.
For any $X\in\ind\C$, by \cite[Proposition 3.12]{Iyama-Nakaoka-Palu_2024}, the sink morphism $u_X: S_X \to X$ induces an isomorphism 
\begin{align}
    (\C/\rad\C)(-, S_X) \xrightarrow[\sim]{(-, u_X)}  (\rad\C/\rad^2\C)(-, X)= \rad\C(-, X).
\end{align}
This shows that $\rad\C(-, X)$ is a finitely generated semisimple functor.
Note that every finitely generated semisimple functor $F\in\Mod\C$ has finite support and satisfies that 
$F(Y)_{\C(Y, Y)}$ has finite length for all $Y\in\ind\C$.
Together with this observation, the short exact sequence 
\begin{align}
    0 \to \rad\C(-,X) \to \C(-, X)\to (\C/\rad\C)(-, X) \to 0
\end{align}
shows that $\C(-, X)$ has finite support, and $\C(Y,X)_{\C(Y,Y)}$ has finite length for all $Y\in\ind\C$.
Therefore, $\Mod\C$ is locally finite.
\end{proof}

Let $X$ be an object of $\C$.
If the projective functor $\C(-, X)\in\Mod\C$ has finite length, that is, there is a finite filtration $0=F_0 \subset F_1 \subset \cdots  \subset F_n=\C(-, X)$ such that $F_{i+1}/F_i$ is simple for $0\leq i < n$, then we denote its length $n$ by $l(X)$.
We write $l(X)=\infty$ if $\C(-, X)$ has infinite length.
If $\Mod\C$ is locally finite, then $l(X)< \infty$ for all objects $X\in\C$ (\cite[Proposition 1.11]{Auslander_1974_II}).
We now give necessary and sufficient conditions for a Krull-Schmidt category with sink morphisms to have square-zero radical.

\begin{proposition}\label{claim_15}
The following conditions are equivalent for a Krull-Schmidt category $\C$ with sink morphisms.
\begin{enumerate}
\item $\rad^2\C = 0$.
\item Every projective functor in $\mod\C$ has either zero or finitely generated semisimple radical.
\item Every sink morphism $u_Z:S_Z \to Z$ admits a weak kernel $u:  X \to S_Z$ such that $\Im(-,u) = \rad\C(-, S_Z)$, where $(-,u): \C(-, X) \to\C(-, S_Z)$.
\end{enumerate}
\end{proposition}
\begin{proof} 
(1) $\Rightarrow$ (2)
Assume (1) holds.
Since $\C$ is Krull-Schmidt, it suffices to show (2) for indecomposable objects in $\C$. 
However, for each $Z\in\ind\C$ with $\rad\C(-, Z)\neq 0$, we know from the proof of \cref{claim_29} that $\rad\C(-, Z)$ is finitely generated semisimple.

(2) $\Rightarrow$ (3)
Assume (2) holds.
It suffices to consider $Z\in\ind\C$ with $\rad\C(-, Z)\neq0$.
For such $Z$, taking a projective cover $\C(-, X) \to \rad\C(-, S_Z)$ of the semisimple functor $\rad\C(-, S_Z)$ and applying Yoneda's lemma, we obtain the desired weak kernel $u: X\to S_Z$ in $\C$ such that $\Im(-,u) = \rad\C(-, S_Z)$, where $(-,u): \C(-, X) \to\C(-, S_Z)$.

(3) $\Rightarrow$ (1)
Assume (3) holds.
We show that $\rad^2\C(-, Z) = 0$ for all $Z \in \ind\C$. 
By assumption, there exists a sequence $X\xrightarrow{u}S_Z\xrightarrow{u_Z}Z$, where $(-,u): \C(-, X) \to\C(-, Y)$ satisfies $\Im(-,u) = \rad\C(-, S_Z)$.
Consider the following induced exact sequence of functors:
\begin{align}
    \C(-, X) \xrightarrow{(-, u)}
    \C(-, S_Z) \xrightarrow{(-, u_Z)}
    \C(-, Z) \to
    (\C/\rad\C)(-, Z) \to 0.
\end{align}
Then we have 
\begin{align}\label{eq_1}
    \rad\C(-, Z) = \Im(-, u_Z)
    &\xrightarrow{\,\sim\,} \Coker(-, u) = (\C/\rad\C)(-, S_Z) \\
    &\xrightarrow[\sim]{\,(-, u_Z)\,} (\rad\C/\rad^2\C)(-, Z). 
\end{align}
A direct calculation shows that the composition of these isomorphisms is the canonical projection from $\rad\C(-, Z)$ to the quotient $(\rad\C/\rad^2\C)(-, Z)$.
This concludes that $\rad^2\C(-, Z) = 0$.
\end{proof}

Next, we characterize radical-square-zero Krull-Schmidt extriangulated categories with sink morphisms.

\begin{proposition}\label{claim_39}
The following conditions are equivalent for a Krull-Schmidt extriangulated
category $(\C, \E, \s)$ with sink morphisms.
\begin{enumerate}
\item $\rad^2\C = 0$.
\item The following two conditions hold.
\begin{enumerate}
    \item Every indecomposable projective functor $\C(-, Z)\in\mod\C$ with $Z\in\ind\C$ non-projective is of length at most two.
    \item Every indecomposable projective functor $\C(-, P)\in\mod\C$ with $P\in\ind\Proj(\C)$ satisfies condition (2) in \cref{claim_15}.
\end{enumerate}
\item The following two conditions hold.
\begin{enumerate}
    \item Every Auslander-Reiten $\E$-triangle $\tri{\tau Z}{}{S_Z}{}{Z}{}$ satisfies $S_Z=0$ or $\rad\C(-, S_Z)\cong (\C/\rad\C)(-,\tau Z)$.
    \item Every indecomposable projective object $P\in\ind\Proj(\C)$ satisfies condition (3) in \cref{claim_15}.
\end{enumerate}
\end{enumerate}
\end{proposition}
\begin{proof} 
(1) $\Rightarrow$ (3)
Assume (1) holds.
By \cref{claim_15}, it suffices to show that (3)(a) holds.
Let $\tri{\tau Z}{u^{\tau Z}}{S_Z}{u_Z}{Z}{}$ be an Auslander-Reiten $\E$-triangle with $S_Z\neq 0$.
Since $\rad^2\C = 0$ by assumption, we have 
\begin{align}
    \Im(-,u^{\tau Z}) = \Ker(-, u_Z) = \Im(-,u) = \rad\C(-, S_Z),
\end{align}
where the first equality follows from \cite[Proposition 3.3]{Nakaoka-Palu_2019_extriangulated} and $u:X\to S_Z$ is the weak kernel appearing in \cref{claim_15}(3).
Since the projective functor $\C(-, \tau Z)$ is indecomposable, the epimorphism $\C(-, \tau Z) \to\rad\C(-, S_Z)$ appearing in the canonical factorization of $(-,u^{\tau Z}): \C(-, \tau Z) \to\C(-, S_Z)$ is a projective cover of the finitely generated semisimple functor $\rad\C(-, S_Z)$.
This shows that $\rad\C(-, S_Z)$ is simple.
Note that this yields $S_Z \in \ind\C$. 
Since $0\neq u^{\tau Z} \in \rad\C(\tau Z, S_Z)$, we infer that $\rad\C(-, S_Z)\cong (\C/\rad\C)(-, \tau Z)$.

(3) $\Rightarrow$ (2)
Assume (3) holds.
We first claim that \cref{claim_15}(3) holds, which implies that (2)(b) holds.
It suffices to show that each non-projective object $Z\in\ind\C$ satisfies \cref{claim_15}(3).
Take an Auslander-Reiten $\E$-triangle $\tri{\tau Z}{u^{\tau Z}}{S_Z}{u_Z}{Z}{}$ ending at $Z$.
We may assume that $S_Z\neq 0$.
Since Yoneda's lemma enables us to obtain
\begin{align}
    \rad(\mod\C)(\C(-, X), \C(-, Y))
    =\{(-, u) :\C(-, X) \to \C(-, Y) \mid u \in \rad\C(X, Y)\}
\end{align}
for any $X, Y\in\C$, \cite[Corollary 3.8]{Krause_2015} implies that the non-zero image $\Im(-, u^{\tau Z})$ coincides with the simple radical $\rad\C(-, S_Z)$.
Therefore, $u^{\tau Z}:\tau Z\to S_Z$ is the desired weak kernel of $u_Z: S_Z \to Z$.
Next, we claim that (3)(a) implies (2)(a).
If $Z\in\ind\C$ is a non-projective object with $\rad\C(-, Z)\neq 0$, then $\rad\C(-, Z)$ is simple.
Indeed, this follows from the following isomorphism of functors
\begin{align}
    \rad\C(-, Z) \cong \Coker(-, u^{\tau Z}) = \C(-, S_Z)/\Im(-, u^{\tau Z}) = (\C/\rad\C)(-, S_Z),
\end{align}
where $S_Z \in \ind\C$ (since $\rad\C(-, S_Z)\cong (\C/\rad\C)(-,\tau Z)$ is simple).

(2) $\Rightarrow$ (1)
Follows from \cref{claim_15}.
\end{proof}

We now characterize radical-square-zero Krull-Schmidt triangulated categories with Auslander-Reiten triangles.
We begin by recalling some necessary notions.

Let $Q=(Q_0, Q_1, s, t)$ be a quiver, where $Q_0$ and $Q_1$ are the sets of vertices and arrows, respectively, and $s, t: Q_1 \to Q_0$ are maps.
For each $a\in Q_1$, we call $s(a)$ and $t(a)$ its {\it source} and {\it target}, respectively.
For each $x \in Q_0$, denote by $x^-$ (resp.~$x^+$) the set of vertices $y\in Q_0$ such that there is an arrow $y\to x$ (resp.~$x \to y$) in $Q$.
A {\it  valued quiver} is a triple $(Q, d, d^\prime)$ of a quiver $Q$ with no multiple arrows and maps $d, d^\prime: Q_1 \to \Z_{\geq 1}$.
We associate to each arrow $x \to y$ in $Q$ the pair $(d_{x, y}, d_{x, y}^\prime)$, called the {\it valuation}.
Note that valued quivers considered in this paper may have loops, whereas those in \cite{Happel-Preiser-Ringel_1980} are not allowed to have loops.
A {\it stable valued translation quiver} is a quadruple $(Q, d, d^\prime, \tau)$, where $(Q, d, d^\prime)$ is a valued quiver, and $\tau: Q_0\to Q_0$ is a bijection, called the {\it translation}, such that $(\tau y)^+ = y^-$ and $(d_{\tau y, x}, d_{\tau y, x}^\prime) = (d_{x, y}^\prime, d_{x, y})$ for all $y \in Q_0$ and $x\in y^-$.
For a group $G$ of automorphisms of $Q$, we define the orbit quiver $Q/G=(Q_0/G,\,  Q_1/G,\,  \overline{s},\,  \overline{t})$, where the maps $\overline{s},\overline{t}:Q_1/G\to Q_0/G$ are defined by $\overline{s}(G\alpha)=Gs(\alpha)$ and $\overline{t}(G\alpha)=Gt(\alpha)$ for $\alpha\in Q_1$.

Recall from \cite{Happel_book} that the {\it Auslander-Reiten quiver} $\G(\T)$ of a Krull-Schmidt triangulated category $\T$ with Auslander-Reiten triangles is the valued quiver defined as follows. 
The vertex set $\G(\T)_0$ is given by $\ind\T$, and we draw a valued arrow $X \xrightarrow{(d_{X,Y},\, d_{X,Y}^\prime)} Y$ for any $X, Y \in \G(\T)_0$ with $d_{XY}\neq 0$, where $d_{X, Y}$ and $d_{X, Y}^\prime$ are defined in \eqref{eq_2} and \eqref{eq_3}.
Moreover, there is a bijection $\tau: \G(\T)_0 \to \G(\T)_0$, called the {\it Auslander-Reiten translation}, such that $X=\tau Z$ if and only if $\T$ has an Auslander-Reiten triangle $X\to Y\to Z \to \Sigma X$.
If $d_{XY}\neq 0$, then \cite[Proposition 3.12]{Iyama-Nakaoka-Palu_2024} implies that there are Auslander-Reiten triangles 
\begin{align}
    \tau Y \to  M \oplus X^{d_{XY}}\to Y \to \Sigma \tau Y
    \quad\text{and}\quad
    X \to Y^{d_{XY}^\prime}\oplus N \to \tau^{-1} X \to \Sigma X.
\end{align}
In particular, $d_{XY} = d_{\tau Y, X}^\prime$.
We say that {\it the Auslander-Reiten quiver of $\T$ is a stable valued translation quiver} if the Auslander-Reiten quiver $(\G(\T), d, d^\prime)$ of $\T$ and its Auslander-Reiten translation $\tau$ form a stable valued translation quiver $(\G(\T), d, d^\prime, \tau)$.

Krause \cite[Lemma A.12]{Krause_2012} proved that if $\T$ is an $R$-linear Hom-finite Krull-Schmidt triangulated category with Auslander-Reiten triangles, then $d_{\tau Y, X} = d_{X, Y}^\prime$ for all $X, Y\in\ind\T$.
This shows that the Auslander-Reiten quiver of such a $\T$ is a stable valued translation quiver.
Moreover, extending a result of Xiao and Zhu \cite[Theorem 2.3.5]{Xiao-Zhu_2005}, Krause proved that every $R$-linear Hom-finite locally finite triangulated category $\T$ has Auslander-Reiten triangles (\cite[Corollary 3.2]{Krause_2012}), and each connected component of $\G(\T)$ is isomorphic to $\Z\Delta/G$ for some Dynkin diagram $\Delta$ and some weakly admissible group $G$ of automorphisms of the translation quiver $\Z\Delta$ (\cite[Theorem 3.5]{Krause_2012}). 
Recall that the following diagrams are the Dynkin diagrams \cite{Happel-Preiser-Ringel_1980}, where unlabelled edges have trivial valuation $(1,1)$:
\[\begin{tikzpicture}[scale=1,inner sep=-.3,outer sep=-.3]
\node (a) at (-.7,0) {$A_n$};
\node (a1) at (0,0) {$\circ$};
\node (a2) at (1,0) {$\circ$};
\node (a3) at (2,0) {$\circ$};
\node (a4) at (2.7,0) {$\cdots$};
\node (a5) at (3.4,0) {$\circ$};
\node (a6) at (4.4,0) {$\circ$};
\node (a7) at (5.5,0) {$\circ$};
\draw[-] (a1) -- (a2) ;
\draw[-] (a2) -- (a3) ;
\draw[-] (a5) -- (a6) ;
\draw[-] (a6) -- (a7) ;
\node (b) at (-.7,-1.) {$B_n$};
\node (b1) at (0,-1.) {$\circ$};
\node (b2) at (1,-1.) {$\circ$};
\node (b3) at (2,-1.) {$\circ$};
\node (b4) at (2.7,-1.) {$\cdots$};
\node (b5) at (3.4,-1.) {$\circ$};
\node (b6) at (4.4,-1.) {$\circ$};
\node (b7) at (5.5,-1.) {$\circ$};
\node at (.5,-.7) {\footnotesize $(1,2)$};
\draw[-] (b1) -- (b2) ;
\draw[-] (b2) -- (b3) ;
\draw[-] (b5) -- (b6) ;
\draw[-] (b6) -- (b7) ;
\node (c) at (-.7,-2.) {$C_n$};
\node (c1) at (0,-2.) {$\circ$};
\node (c2) at (1,-2.) {$\circ$};
\node (c3) at (2,-2.) {$\circ$};
\node (c4) at (2.7,-2.) {$\cdots$};
\node (c5) at (3.4,-2.) {$\circ$};
\node (c6) at (4.4,-2.) {$\circ$};
\node (c7) at (5.5,-2.) {$\circ$};
\node at (.5,-1.7) {\footnotesize $(2,1)$};
\draw[-] (c1) -- (c2) ;
\draw[-] (c2) -- (c3) ;
\draw[-] (c5) -- (c6) ;
\draw[-] (c6) -- (c7) ;
\node (d) at (-.7,-3.1) {$D_n$};
\node (d0) at (0,-2.6) {$\circ$};
\node (d1) at (0,-3.6) {$\circ$};
\node (d2) at (1,-3.1) {$\circ$};
\node (d3) at (2,-3.1) {$\circ$};
\node (d4) at (2.7,-3.1) {$\cdots$};
\node (d5) at (3.4,-3.1) {$\circ$};
\node (d6) at (4.4,-3.1) {$\circ$};
\node (d7) at (5.5,-3.1) {$\circ$};
\draw[-] (d0) -- (d2) ;
\draw[-] (d1) -- (d2) ;
\draw[-] (d2) -- (d3) ;
\draw[-] (d5) -- (d6) ;
\draw[-] (d6) -- (d7) ;
\node (e6) at (6.6,-0) {$E_6$};
\node (e60) at (9.3,.6) {$\circ$};
\node (e61) at (7.3,-0) {$\circ$};
\node (e62) at (8.3,-0) {$\circ$};
\node (e63) at (9.3,-0) {$\circ$};
\node (e65) at (10.3,-0) {$\circ$};
\node (e66) at (11.3,-0) {$\circ$};
\draw[-] (e60) -- (e63) ;
\draw[-] (e61) -- (e62) ;
\draw[-] (e62) -- (e63) ;
\draw[-] (e63) -- (e65) ;
\draw[-] (e65) -- (e66) ;
\node (e7) at (6.6,-1.) {$E_7$};
\node (e70) at (9.3,-.4) {$\circ$};
\node (e71) at (7.3,-1.) {$\circ$};
\node (e72) at (8.3,-1.) {$\circ$};
\node (e73) at (9.3,-1.) {$\circ$};
\node (e75) at (10.3,-1.) {$\circ$};
\node (e76) at (11.3,-1.) {$\circ$};
\node (e77) at (12.3,-1.) {$\circ$};
\draw[-] (e70) -- (e73) ;
\draw[-] (e71) -- (e72) ;
\draw[-] (e72) -- (e73) ;
\draw[-] (e73) -- (e75) ;
\draw[-] (e75) -- (e76) ;
\draw[-] (e76) -- (e77) ;
\node (e8) at (6.6,-2.) {$E_8$};
\node (e80) at (9.3,-1.4) {$\circ$};
\node (e81) at (7.3,-2.) {$\circ$};
\node (e82) at (8.3,-2.) {$\circ$};
\node (e83) at (9.3,-2.) {$\circ$};
\node (e85) at (10.3,-2.) {$\circ$};
\node (e86) at (11.3,-2.) {$\circ$};
\node (e87) at (12.3,-2.) {$\circ$};
\node (e88) at (13.3,-2.) {$\circ$};
\draw[-] (e80) -- (e83) ;
\draw[-] (e81) -- (e82) ;
\draw[-] (e82) -- (e83) ;
\draw[-] (e83) -- (e85) ;
\draw[-] (e85) -- (e86) ;
\draw[-] (e86) -- (e87) ;
\draw[-] (e87) -- (e88) ;
\node (f) at (6.6,-3.1) {$F_4$};
\node (f1) at (7.3,-3.1) {$\circ$};
\node (f2) at (8.3,-3.1) {$\circ$};
\node (f3) at (9.3,-3.1) {$\circ$};
\node (f4) at (10.3,-3.1) {$\circ$};
\node at (8.8,-2.8) {\footnotesize $(1,2)$};
\draw[-] (f1) -- (f2) ;
\draw[-] (f2) -- (f3) ;
\draw[-] (f3) -- (f4) ;
\node (g) at (6.6,-4.1) {$G_2$};
\node (g1) at (7.3,-4.1) {$\circ$};
\node (g2) at (8.3,-4.1) {$\circ$};
\node at (7.8,-3.8) {\footnotesize $(1,2)$};
\draw[-] (g1) -- (g2) ;
\end{tikzpicture}
\]
Recall also that a group $G$ of automorphisms of a quiver $Q$ is called {\it weakly admissible} \cite[Section 4]{Dieterich_2006} if $x^+ \cap (gx)^+ = \varnothing$ for all $g \in G\setminus\{1\}$ and $x \in Q_0$.
Amiot \cite[Theorem 2.1]{Amiot_2007} proved that if $\Delta$ is a Dynkin diagram of type A, D or E, then any non-trivial weakly admissible group $G$ of automorphisms of $\Z\Delta$ is isomorphic to $\Z$.
She also gave a generator of $G$.

\cref{claim_3} below may be regarded as a result similar to \cite[Theorem 3.5]{Krause_2012} for (not necessarily $R$-linear Hom-finite) radical-square-zero Krull-Schmidt triangulated categories with Auslander-Reiten triangles.
Let us begin with a series of lemmas.

\begin{lemma}\label{claim_26}
Assume that $\rad^2\T=0$ and $\T$ has right Auslander-Reiten triangles.
\begin{enumerate}
\item Every irreducible morphism between indecomposable objects is both a sink and a source morphism. 
\item  If $X \xrightarrow{u} Y \xrightarrow{v} Z \xrightarrow{w} \Sigma X$ is an Auslander-Reiten triangle with $Y\in\ind\T$, then the triangles 
\begin{align}
   Y \xrightarrow{v} Z \xrightarrow{w} \Sigma X \xrightarrow{-\Sigma u} \Sigma Y
   \quad\text{and}\quad
   \Sigma^{-1} Z \xrightarrow{-\Sigma^{-1}w} X \xrightarrow{u} Y \xrightarrow{v} Z
\end{align}
are also Auslander-Reiten triangles.
\item  $\T$ has left Auslander-Reiten triangles.
\end{enumerate}
\end{lemma}
\begin{proof}
(1) Let $f:X\to Y$ be an irreducible morphism between indecomposables.
Then the sink morphism $u_Y: S_Y \to Y$ contains $f$ as a direct summand.
Since $S_Y$ is indecomposable by \cref{claim_39}, it follows that $f$ is isomorphic to $u_Y$.
The dual argument shows that $f:X\to Y$ is also isomorphic to the source morphism $u^X: X \to S^X$.

(2) Let $X \xrightarrow{u} Y \xrightarrow{v} Z \xrightarrow{w} \Sigma X$ be an Auslander-Reiten triangle with $Y\in\ind\T$.
Since $u=u^X=u_Y$ and $v=u_Z=u^Y$ by (1), we obtain the desired result (cf.~\cite[Theorem 2.9]{Iyama-Nakaoka-Palu_2024}).

(3) Follows from (2) 
\end{proof}

\begin{lemma}\label{claim_31}
Assume that $\rad^2\T=0$ and $\T$ has Auslander-Reiten triangles.  
\begin{enumerate}
\item   For any $X, Y\in\ind\T$, we have 
\begin{align}
d_{X, Y} = 
\begin{cases}
    \delta_{X, S_Y} & \text{ if $S_Y\neq 0$},\\
    0 & \text{ if $S_Y =0$}
\end{cases}
    \quad\text{and}\quad
d_{X, Y}^\prime = 
\begin{cases}
    \delta_{Y, S^X} & \text{ if $S^X\neq 0$},\\
    0 & \text{ if $S^X =0$},
\end{cases}
\end{align}
where $ \delta_{A, B} = 1$ if $A \cong B$ and otherwise $\delta_{A, B} = 0$. 
\item 
For any $X, Y\in\ind\T$, we have $d_{\tau Y, X} = d_{X, Y}^\prime$.
In particular, the Auslander-Reiten quiver of $\T$ is a stable valued translation quiver.
\end{enumerate}
\end{lemma}
\begin{proof} 
(1) Let $X, Y\in\ind\T$ be arbitrary.
For the first equality, it is clear that $d_{X, Y}=0$ if $S_Y=0$.
Assume $S_Y\neq 0$.
Since the Auslander-Reiten triangle $\tau Y \to S_Y \to Y \to \Sigma \tau Y$ ending at $Y$ satisfies $S_Y\in\ind\T$ by \cref{claim_39}, it follows that $d_{X, Y} = \delta_{X, S_Y}$ as desired.
Dually for the second equality.

(2) For any $X, Y\in\ind\T$, one can use the Auslander-Reiten triangle $\tau Y \xrightarrow{u^{\tau Y}=u_{S_Y}} S_Y \xrightarrow{u_Y=u^{S_Y}} Y \to \Sigma \tau Y$ ending at $Y$ to deduce that $\tau Y \cong S_X$ if and only if $Y \cong S^X$.
Then (1) yields the desired equality $d_{\tau Y, X} = d_{X, Y}^\prime$.
\end{proof}

Let $(Q, d, d^\prime, \tau)$ be a valued translation quiver.
Following \cite{Happel-Preiser-Ringel_1980}, we say that 
a function $l: Q_0 \to \Z_{\geq 1}$ is {\it subadditive} if $l(x) + l(\tau x) \geq \sum_{y\in x^-} l(y) d_{y, x}^\prime$ for all $x\in Q_0^\prime$.
If the equality holds for all $x\in Q_0^\prime$, we then call $l$ {\it additive}.

\begin{lemma}\label{claim_32}
Let $\T$ be as in \cref{claim_31}.
\begin{enumerate}
\item  $l(X) \leq 2$ for all $X\in\ind\T$. 
The equality holds if and only if $S_X\neq 0$.
\item 
The function $d: \G(\T)_0 \to \Z_{\geq 1}$ given by $Z\mapsto l(Z)$ has the property that for each $Z\in\G(\T)_0$, we have 
\begin{align}
    2 d(Z) = d(Z) + d(\tau Z) = 2+ \sum_{Y\in Z^-} d_{Y, Z} d(Y).
\end{align}
\end{enumerate}
\end{lemma}
\begin{proof} 
(1) For any $X\in\ind\T$, we have $l(X) <\infty$ since $\Mod\T$ is locally finite by \cref{claim_29}.
Since $S_X$ is zero or indecomposable by \cref{claim_39}, it follows that $\rad\T(-, X)=( \rad\T/\rad^2\T)(-, X)\cong(\T/\rad\T)(-, S_X)$ is zero or simple.
Thus the short exact sequence
\begin{align}
    0 \to \rad\T(-,X) \to \T(-, X)\to (\T/\rad\T)(-, X) \to 0
\end{align}
concludes that $l(X) \leq 2$, and the equality holds if and only if $S_X\neq 0$.

(2) For any $Z\in\ind\T$, the second equality follows from the exact sequence 
\begin{align}\label{eq_4}
    0 \to (\T/\rad\T)(-, \Sigma^{-1}Z) \to &\T(-,\tau Z) \xrightarrow{(-, u)}& \\
    &\T(-, S_Z)\xrightarrow{(-, v)} \T(-,Z) \to (\T/\rad\T)(-, Z) \to 0
\end{align}
induced by the Auslander-Reiten triangle $\tau Z \xrightarrow{u} S_Z \xrightarrow{v} Z \xrightarrow{w} \Sigma\tau Z$ ending at $Z$ (cf.~the proof of \cite[Proposition 3.4]{Krause_2012}).
For the first equality, if $S_Z=0$, then $\T(-, \tau Z)$ and $\T(-, Z)$ are both simple by \eqref{eq_4}, which shows that $d(\tau Z) =1 = d(Z)$.
If $S_Z\neq 0$, then $u=u^{\tau Z}=u_{S_Z}$ and $-\Sigma^{-1}w=u^{\Sigma^{-1}Z}=u_{\tau Z}$ by \cref{claim_26}. 
Since $S_{\tau Z}=\Sigma^{-1}Z\neq0$, it follows from (1) that  $d(\tau Z) = 2 = d(Z)$.
\end{proof}

We are now ready to prove the main result of this section.

\begin{theorem}
\label{claim_3}
The following conditions are equivalent for a Krull-Schmidt triangulated category $\T$ with Auslander-Reiten triangles.
\begin{enumerate}
\item $\rad^2\T = 0$.
\item Every indecomposable projective object in $\mod\T$ is of length at most two.
\item The middle term of every Auslander-Reiten triangle in $\T$ is either zero or indecomposable.
\item Every connected component of the Auslander-Reiten quiver of $\T$ is of the form $\ZA_n/G$ for some $n \in\{1, 2\}$ and some weakly admissible group $G$ of automorphisms of $\ZA_n$.  
\end{enumerate}
\end{theorem}
\begin{proof} 
The equivalences (1) $\Leftrightarrow$ (2) $\Leftrightarrow$ (3) follows from \cref{claim_39} and the fact that $\Proj(\T)=\{0\}$.

(1) $\Rightarrow$ (4)
Assume $\rad^2\T = 0$.
Using \cref{claim_29} and Lemmas \ref{claim_31}(2) and \ref{claim_32}(2) , one can prove, as in the proof of \cite[Theorem 2.3.5]{Xiao-Zhu_2005} (cf.~that of \cite[Theorem 4.1]{Amiot_2007}), that the Auslander-Reiten quiver $\G(\T)$ of $\T$ is a stable valued translation quiver, and each connected component is of the from $\Z\Delta/G$ for some Dynkin diagram $\Delta$ and some weakly admissible group $G$ of automorphisms of $\Z\Delta$.
On the other hand, \cref{claim_31}(1) implies that for each vertex $X \in\G(\T)_0$, there is at most one arrow in $\G(\T)$ starting at $X$ with valuation $(1, 1)$, and there is at most one arrow in $\G(\T)$ ending at $X$ with valuation $(1, 1)$.
This forces $\Delta = A_n$ with $1 \leq n \leq 2$.

(4) $\Rightarrow$ (3)
We first claim (without condition (4)) that an Auslander-Reiten triangle $X \to Y \to Z \to \Sigma X$ in $\T$ satisfies $Y = 0$ if and only if $Z \in\G(\T)_0$ is both a sink and a source in $\G(\T)$.
If $Y=0$, then, as in the proof of \cref{claim_26}, one shows that $\Sigma^{-1}Z \to 0 \to Z \xrightarrow{1} Z$ and $Z \to 0 \to \Sigma Z \xrightarrow{1} \Sigma Z$ are both Auslander-Reiten triangles.
This shows that there is no arrow either starting or ending at $Z \in\G(\T)_0$ in $\G(\T)$.
The converse is an easy consequence of \cite[Proposition 3.12]{Iyama-Nakaoka-Palu_2024}.

Assume now (4) holds.
Let $Z \in\ind\T$ lie in a connected component (isomorphic to) $\ZA_n/G$ for some $n \in\{1, 2\}$ and some weakly admissible group $G \leq \Aut(\ZA_n)$. 
If $n=1$, then $S_Z=0$, since $\ZA_1/G$ is discrete (i.e.~it has no arrow).
Suppose $n=2$.
By \cite[Theorem 2.1]{Amiot_2007}, $G$ is the cyclic group generated by an automorphism of the form $\tau^{-r}$ or $\tau^{-r}\Sigma$ for some $r\geq 0$. 
In particular, $G$ is non-trivial if and only if $\ZA_2/G$ has only  finitely many vertices.
In this case, $\ZA_2/G$ is a cyclic quiver.
Thus, in $\ZA_2/G$, there is exactly one arrow ending at $Z$ with valuation $(1, 1)$.
Together with this observation, \cite[Proposition 3.12]{Iyama-Nakaoka-Palu_2024} concludes that $S_Z \in \ind\T$.
 \end{proof}

\begin{example}
Let $\Delta$ be a Dynkin diagram of type A, D or E. 
Since the Auslander-Reiten quiver of $\D^{\rm b}(\mod k\Delta)$ is of the form $\Z\Delta$ (\cite[I.5.5]{Happel_book}), it follows from \cref{claim_3} that $\rad^2\D^{\rm b}(\mod k\Delta)=0$ if and only if $\Delta
=A_n$ for some $1\leq n \leq 2$.
Moreover, $\rad\D^{\rm b}(\mod k\Delta)=0$ if and only if $\Delta=A_1$ (cf.~\cref{claim_5}).
\end{example}

We say that $\T$ is {\it connected} (or a {\it block}) if it cannot be decomposed into a product of two non-zero triangulated categories. 
A triangulated subcategory of $\T$ is called {\it thick} if it is closed under direct summands.
Every thick subcategory of $\T$ is a Krull-Schmidt category.
For an object $T \in \T$, we denote by $\thick(T)=\thick_\T(T)$ the smallest thick subcategory of $\T$ containing $T$.
Following \cite[4.1]{Chen-Ringel_2018}, we call a sequence $X_0, X_1, \ldots,  X_n$ of indecomposable objects in $\T$ such that $\T(X_{i-1},X_i) \neq 0$ or $X_i = \Sigma X_{i-1}$ for $1 \leq i \leq n$ a {\it path from $X_0$ to $X_n$}.
A subset $\U \subseteq \ind\T$ is called {\it path-closed} if for each path from $X$ to $Y$, $X\in\U$ if and only if $Y\in\U$.
Observe that $\U$ is path-closed if and only if $\U$ is closed under the shift functors $\Sigma$ and $\Sigma^{-1}$, and for any $X\in\U$, every $Y\in\ind\T$ with $\T(X, Y)\neq 0$ or $\T(Y, X)\neq 0$ lies in $\U$.
We say that $\U$ is {\it path-connected} if any two objects in $\U$ are connected by a sequence of paths and inverse paths, more precisely, for any $X, Y\in\U$, there exists a sequence $X = X_0, X_1, \ldots, X_t = Y$ of indecomposable objects in $\T$ such that for $1 \leq i \leq t$, at least one of the three conditions $\T(X_{i-1}, X_i)\neq 0$, $\T(X_i, X_{i-1})\neq 0$ and $X_i = \Sigma^sX_{i-1}$ for some $s \in \Z$, is satisfied.

\begin{proposition}\label{claim_40}
Let $\T$ be a radical-square-zero Krull-Schmidt triangulated category with Auslander-Reiten triangles.
\begin{enumerate}
\item For each path-closed subset $\U$ of $\ind\T$, its additive closure $\T_\U=\add\U$ in $\T$ is a radical-square-zero connected thick subcategory of $\T$ with Auslander-Reiten triangles such that  
\begin{align}
    \T_\U = \thick_\T(S_X\oplus X)
\end{align}
for any $X\in\U$. 
\item For any two distinct path-closed subsets $\U, \U^\prime$ of $\ind\T$, we have 
\begin{align}
    \T(\T_{\U}, \T_{\U^\prime}) =0= \T(\T_{\U^\prime}, \T_{\U}).
\end{align}
In particular, if $\ind\T$ has a finite number of path-closed subsets, then there is a decomposition $\T\,=\,\prod_{\U}\T_\U$ into a finite product of connected thick subcategories, where $\U$ runs over the path-closed subsets of $\ind\T$.
\end{enumerate}
\end{proposition}
\begin{proof}
Since $\rad^2\T=0$, \cref{claim_26} implies that if $A, B \in\ind\T$, then $\rad\T(A, B) \neq 0$ if and only if $A\cong S_B$ if and only if  $B\cong S^A$.
In this case, any non-zero non-isomorphism $f\in \rad\T(A, B)$ is both a sink and a source morphism.

(1) Let $\U$ be a path-closed subset of $\ind\T$.
Then the additive closure $\T_\U=\add\U$ is a thick subcategory of $\T$ by \cite[Lemma 4.1]{Chen-Ringel_2018}.
Moreover, since $\ind\T_\U = \U$ is path-connected, \cite[Proposition 4.2]{Chen-Ringel_2018} implies that $\T_\G$ is connected.
Further, the fact that every Auslander-Reiten triangle in $\T$ whose terms are all in $\U$ becomes an Auslander-Reiten triangle in $\T_\U$ implies that $\T_\U$ has Auslander-Reiten triangles and satisfies $\rad^2\T_\U=0$.
Now, for each $X\in\U$, consider the thick closure $\thick_\T(S_X\oplus X)$ of $S_X\oplus X$ in $\T$.
Then the observation at the beginning of the proof shows that $\ind\thick(S_X\oplus X)=\U$, so that $\thick_\T(S_X\oplus X)=\T_\U$.

(2) Again, the observation at the beginning shows that  $\T(\T_\U, \T_{\U^\prime}) =0= \T(\T_{\U^\prime}, \T_\U)$ for any two distinct path-closed subsets $\U, \U^\prime$ of $\ind\T$ (otherwise, we have $\U=\U^\prime$, a contradiction).
This completes the proof.
\end{proof}

A triangulated category is called {\it algebraic} if it is triangle equivalent to the stable category of a Frobenius exact category.

We say that a $\mathbbm{k}$-linear Hom-finite Krull-Schmidt triangulated category $\T$ {\it has a Serre duality} if there exists a triangle auto-equivalence $\nu: \T \to \T$, called a {\it Serre functor}, together with binatural isomorphisms $\T(X, Y) \cong D\T(Y, \nu X)$ for all $X, Y\in \T$, where $D=\Hom_\mathbbm{k}(-, E)$ denotes the Matlis duality, where $E$ is an injective envelope of $\mathbbm{k}/\rad \mathbbm{k}$.
A Serre functor is unique up to natural isomorphism if it exists.
Moreover, an argument as in \cite[Section 3]{Iyama-Nakaoka-Palu_2024} (cf.~\cite[Section 6.3]{Krause_book}) shows that $\T$ has a Serre duality if and only if $\T$ has Auslander-Reiten triangles.
In this case, $\nu X = \Sigma \tau X$ for all $X\in\ind\T$.
It is known that the bounded derived category $\D^{\rm b}(\mod \L)$ of a finite dimensional $K$-algebra $\L$ with $\gldim\L<\infty$ has a Serre duality (\cite[3.4]{Happel_1991}).
Let $l\in\Z$.
We say that $\T$ is {\it $l$-Calabi-Yau} \cite{Keller_2005} if $\Sigma^l$ defines a Serre duality.
In the rest of this paper, $\nu$  always denotes a Serre functor of a triangulated category with Serre duality.

For a triangle autoequivalence $\Phi$ of $\D^{\rm b}(\mod K\Delta)$
with $\Delta$ a Dynkin diagram of type ADE, we denote by 
$\D^{\rm b}(\mod K\Delta)/\Phi$ the triangulated orbit category of $\D^{\rm b}(\mod K\Delta)$ by the cyclic group generated by $\Phi$ in the sense of Keller \cite{Keller_2005}.
Thanks to Amiot \cite{Amiot_2007}, we obtain the following.

\begin{proposition}\label{claim_36}
Keep the notation in \cref{claim_40}.
Assume further that $\T$ is $K$-linear Hom-finite and algebraic.
Then, for each path-closed subset $\U$ of $\ind\T$, there is a triangle equivalence  
\begin{align}
    \T_\U \cong  \D^{\rm b}(\mod KA_n)/\Phi
\end{align}
for some $n \in\{1, 2\}$ and some triangle auto-equivalence $\Phi$ of $\D^{\rm b}(\mod KA_n)$.
\end{proposition}
\begin{proof}
Let $\U$ be a path-closed subset of $\ind\T$.
By the additional assumption and \cref{claim_40}, $\T_\U$ is a radical-square-zero, $K$-linear, Hom-finite, algebraic, Krull-Schmidt, connected triangulated category with Auslander-Reiten triangles.
Then \cref{claim_26} implies that for $A, B \in\ind\T_\U$, we have
\begin{align}\label{eq_12}
    \T_\U(A, B) = 
    \begin{cases}
        K 1_A \oplus \rad\T_\U(A, B) & \text{if } A \cong B,\\
        \rad\T_\U(A, B) & \text{otherwise},
    \end{cases}
\end{align}
where 
\begin{align}
    \rad\T_\U(A, B) =
    \begin{cases}
        K u_B= K u^A & \mbox{if $A\cong S_B$, equivalently, $B\cong S^A$},\\
        0 & \mbox{otherwise.}
    \end{cases}
\end{align}
If $\rad\T_\U=0$, then \cref{claim_37} yields a triangle equivalence $\T_\U\cong \D^{\rm b}(\mod K)/\Sigma^n$ for some $n\in\Z$.
If $\rad\T_\U\neq 0$, then, since $\T_\U$ is connected, we have $\G(\T_\U)\cong \ZA_2/G$ for some weakly admissible automorphism group $G$ of $\ZA_2$ by \cref{claim_3}.
If $\ind\T_\U$ is infinite (i.e.~$G$ is trivial), then \cite[Theorem 5.1]{Krause_2012} implies that $\T_\G\cong \D^{\rm b}(\mod KA_2)$.
Assume now that $\ind\T_\U$ is finite.
Then $\T_\U$ is standard, that is, there is an equivalence $\ind\T_\U \cong K\G(\T_\U)$ of $K$-linear categories, where $K\G(\T_\U)$ denotes the mesh category of the Auslander-Reiten quiver $\G(\T_\U)$ of $\T_\U$ (see \cite[5.6]{Happel_book} for its definition).
Indeed, together with the fact that $\G(\T_\U)$ is a cyclic quiver (see the proof of \cref{claim_3}), \eqref{eq_12} shows that $K\G(\T_\U)$ is $K$-linearly equivalent to $\ind\T_\U$.
Therefore, \cite[Theorem 7.2]{Amiot_2007} yields a triangle equivalence $\T_\U\cong \D^{\rm b}(\mod KA_2)/\Phi$ for some (non-trivial) triangle auto-equivalence $\Phi$ of $\D^{\rm b}(\mod KA_2)$.
\end{proof}

For a finite dimensional $K$-algebra $\L$, we denote by $\Out(\D^{\rm b}(\mod \L))$ the group of triangle auto-equivalences of $\D^{\rm b}(\mod \L)$ up to natural isomorphisms.
By \cite[Corollary 1.9 and Theorem 4.1]{Miyachi-Yekutieli_2001}, we have 
\begin{align}
    \Out(\D^{\rm b}(\mod KA_n)) =\langle\tau, \Sigma\mid \tau^{-n-1}= \Sigma^2\rangle.
\end{align}
Together with the formula $\nu=\Sigma\tau$, a direct calculation shows that 
\begin{align}
\Out(\D^{\rm b}(\mod KA_n)) =
\begin{cases}
    \langle\Sigma\rangle \text{ with } \Sigma^{-1}=\tau & \text{ if $n=1$},\\[0mm]
    \langle\nu\rangle \text{ with $\tau=\nu^{-2}$ and $\Sigma=\nu^3$}  & \text{ if $n=2$}.
\end{cases}
\end{align}

Observe that the shift functor $\Sigma$ (resp.~the Auslander-Reiten translation $\tau$, the Serre functor $\nu$) for $\D^{\rm b}(\mod KA_n)$ induces the shift functor $\underline{\Sigma}$ (resp.~the Auslander-Reiten translation $\underline{\tau}$, the Serre functor $\underline{\nu}$) for $\D^{\rm b}(\mod KA_n)/\Phi$.
Note that $KA_1=K$.

\begin{proposition}\label{claim_42}
For each $i\in\Z$, the following statements hold. 
\begin{enumerate}
\item $\D^{\rm b}(\mod K)/\Sigma^i$ is $0$-Calabi-Yau.
\item $\D^{\rm b}(\mod KA_2)/\nu^i$ is $l$-Calabi-Yau if and only if $3l\equiv 1\ (\mathrm{mod}\ i)$.
\end{enumerate}
\end{proposition}
\begin{proof} 
Keep the notation preceding the lemma.

(1) Follows from the equalities $\underline{\nu}= \underline{\Sigma}\,\underline{\tau}=\underline{\Sigma}\,\underline{\Sigma}^{-1}=1_{\D^{\rm b}(\mod K)/\Sigma^i}$.

(2) For the \lq\lq if\rq\rq~part, assume $3l-1=mi$ for some $m\in\Z$.
Since  $\underline{\nu}^3=\underline{\Sigma}$ and $\underline{\nu}^i\cong 1_{\D^{\rm b}(\mod KA_2)/\nu^i}$, we have $\underline{\nu} = \underline{\nu}^{3l-mi} = \underline{\nu}^{3l}\,\underline{\nu}^{-mi} \cong \underline{\Sigma}^l$, which means that $\D^{\rm b}(\mod KA_2)/\nu^i$ is $l$-Calabi-Yau.
Conversely, if $\underline{\nu} \cong \underline{\Sigma}^l$ for some $l\in\Z$, then $\underline{\nu}^{3l-1} = \underline{\nu}^{3l}\,\underline{\nu}^{-1} \cong \underline{\Sigma}^l\,\underline{\Sigma}^{-l}\cong 1_{\D^{\rm b}(\mod KA_2)/\nu^i}$.
Since $\underline{\nu}^i\cong 1_{\D^{\rm b}(\mod KA_2)/\nu^i}$, it follows that for any $X\in\ind\D^{\rm b}(\mod KA_2)$, there exists an integer $m=m_X\in\Z$ such that $\nu^{3l-1} X \cong \nu^{im}X$ in $\D^{\rm b}(\mod KA_2)$.
Hence $3l-1=im$.
\end{proof}

\begin{remark}\label{remark_3}
Assume that $\T$ is a $K$-linear Hom-finite Krull-Schmidt radical-square-zero triangulated category with Auslander-Reiten triangles such that $\ind\T$  has only finitely many path-closed subsets.
Then \cref{claim_40,claim_36} imply that 
\begin{align}    
\T\,=\,\prod_{\U}\T_\U \,\cong\, \prod_{\U} \D^{\rm b}(\mod KA_{n_\U})/\Phi_\U
\end{align}
as triangulated categories, where $n_\U\in\{1,2\}$ and $\Phi_\U\in \Out(\D^{\rm b}(\mod KA_{n_\U}))$.
One can check whether $\T$ is Calabi-Yau or not by using \cref{claim_42} (cf.~\cref{example_8}).
\end{remark}



\section{G-radical-square-zero algebras} 
\label{G-radical_square_zero_algebras}

Throughout this section, assume that $\L$ is an artin $\mathbbm{k}$-algebra. 
This section is devoted to studying artin algebras $\L$ for which the stable category of Gorenstein-projective modules has square-zero radical.
We begin by recalling the definition of Gorenstein-projective modules and some necessary facts.

A complex  $P^{\bullet} : \cdots \to P^{i-1} \xrightarrow{d^{i-1}} P^{i} \xrightarrow{d^{i}} P^{i+1} \to \cdots$ of projective $\L$-modules is called {\it totally acyclic} \cite{Avramov-Martsinkovsky_2002} if both $P^{\bullet}$ and $\Hom_{\L}(P^\bullet, \L)$ are acyclic. 
A $\L$-module $M$ is called  {\it Gorenstein-projective} \cite{Enochs-Jenda_1995_MathZ} if there exists a totally acyclic complex $P^{\bullet}$ such that $\Ker d^{0}=M$.
Gorenstein-projective modules are  called {\it modules of Gorenstein dimension zero} \cite{Auslander-Bridger_1969}, {\it totally reflexive modules} \cite{Avramov-Martsinkovsky_2002} and {\it maximal Cohen-Macaulay modules} \cite{Buchweitz_1986} in the literature.
We denote by $\Gproj\L$ the subcategory of Gorenstein-projective $\L$-modules.
Since projective modules are Gorenstein-projective, we have $\proj\L \subseteq \Gproj \L \subseteq \mod \L$.
Note that $\Gproj \L = \mod \L$ if and only if $\L$ is self-injective (i.e.~$\L_\L$ is injective).
We say that $\L$ is {\it CM-free} if $\proj\L = \Gproj \L$.
For example, algebras of finite global dimension are CM-free.
We say that $\L$ is {\it CM-finite} if there are only finitely many pairwise non-isomorphic indecomposable Gorenstein-projective $\L$-modules.
Clearly, CM-free algebras are CM-finite.

Let $\mathcal{P}$ be the ideal of $\mod\L$ such that for any $M, N \in \mod\L$,  $\mathcal{P}(M, N)$ is the space of morphisms from $M$ to $N$ that factor through a projective $\L$-module.
The {\it stable category} $\smod\L$ of $\mod\L$ is the ideal quotient of $\mod\L$ by the ideal $\mathcal{P}$. 
Let $\sGproj \L$ be the subcategory of $\smod\L$ consisting of Gorenstein-projective $\L$-modules.
It is clear that $\sGproj \L = 0$ if and only if  $\L$ is CM-free.
The category $\Gproj \L$ is known to be a Frobenius exact category whose projective objects are precisely the projective $\L$-modules.
Therefore, the stable category $\sGproj\L$ carries a structure of a triangulated category; see \cite{Happel_book}.
The functor $(-)^*=\Hom_\L(-, \L)$ induces a duality  $(-)^*: \Gproj\L \to \Gproj\L^\op$ of Frobenius exact categories with quasi-inverse $(-)^*=\Hom_{\L^\op}(-, \L)$ (\cite[Corollary 2.1.6]{X-WChen_2017}).
Hence we obtain a duality $(-)^*: \sGproj\L \to \sGproj\L^\op$ of triangulated categories.
Recall from \cite{Buchweitz_1986,Orlov_2004} that the {\it singularity category} $\D_{\rm sg}(\L)$ of $\L$ is the Verdier quotient of the bounded derived category $\D^{\rm b}(\mod \L)$ of $\mod\L$ by the category $\K^{\rm b}(\proj \L)$ of perfect complexes (i.e.~complexes that are isomorphic in $\D^{\rm b}(\mod\L)$ to bounded complexes of projective $\L$-modules).  
We call $\L$ {\it $d$-Iwanaga-Gorenstein} \cite{Iwanaga_1979} (or simply {\it Iwanaga-Gorenstein}) if both $\injdim_{\L}\L$ and $\injdim \L_{\L}$ are finite and at most $d$. 
In this case, we have $\injdim_{\L}\L=\injdim \L_{\L}$ (\cite[Lemma A]{Zaks_1969}).
By \cite[Theorem 4.4.1]{Buchweitz_1986}, there is a triangle embedding $\Phi: \sGproj \L\hookrightarrow\D_{\rm sg}(\L)$ sending each object $M\in\sGproj\L$ to the stalk complex $M$ concentrated in degree $0$.
Moreover, $\Phi: \sGproj \L\hookrightarrow\D_{\rm sg}(\L)$ is an equivalence if and only if $\L$ is Iwanaga-Gorenstein.

The following class of artin algebras was introduced by Hafezi and Bahlekeh \cite{Hafezi-Bahlekeh_2024}.

\begin{definition}[{\cite[Definition 3.1]{Hafezi-Bahlekeh_2024}}] \label[defi]{def_G-semisimple}
We say that an artin algebra $\L$ is  {\it G-semisimple} if $\mod\sGproj\L$ is a semisimple abelian category.
\end{definition}

\begin{remark} \label{remark_2}
By \cref{claim_6,claim_5}, an artin algebra $\L$ is G-semisimple if and only if $\rad\sGproj\L = 0$ if and only if $\sGproj\L$ is a semisimple triangulated category.
\end{remark}

As examples of G-semisimple algebras, there are quadratic monomial algebras, including gentle algebras.
This is a consequence of \cref{claim_28}, which summarizes known characterizations of G-semisimple monomial algebras established by Chen, Shen and Zhou \cite{Chen-Shen-Zhou_2018} and the authors \cite{Honma-Usui_2025}.

In this section, we extend results of Hafezi and Bahlekeh \cite[Section 3]{Hafezi-Bahlekeh_2024} on G-semisimple algebras to the following class of artin algebras.

\begin{definition}
An artin algebra $\L$ is called {\it G-radical-square-zero} if $\rad^2\sGproj\L$=0.
\end{definition}

G-semisimple algebras are clearly G-radical-square-zero algebras.
Before studying G-radical-square-zero algebras, we focus on artin algebras $\L$ satisfying $\rad^2\Gproj\L = 0$.

For a collection $\X$ of objects in an additive category $\C$, a {\it right $\X$-approximation} of an object $C\in\C$ is a morphism $f : X \to C$ in $\C$ with $X \in \X$ such that for every morphism $f^\prime : X^\prime \to C$ with $X^\prime\in\X$, there is a morphism $g : X^\prime \to X$ such that $f^\prime = fg$.
We simply call a right minimal right $\X$-approximation a {\it minimal right $\X$-approximation}.
We say that $\X$ is {\it contravariantly finite} if any object in $\C$ has a right $\X$-approximation. 
We define a {\it left $\X$-approximation}, a  {\it minimal left $\X$-approximation} and a {\it covariantly finite subcategory} dually. 
We say that $\X$ is {\it functorially finite} if it is contravariantly and covariantly finite.

\begin{proposition}\label{claim_22}
Let $\L$ be an artin algebra.
\begin{enumerate}
\item $\rad^2\Gproj\L=0$ if and only if  $\L$ is a CM-free algebra with radical square zero.
\item $\rad\Gproj\L=0$  if and only if $\L$ is a semisimple algebra.
\end{enumerate}
\end{proposition}
\begin{proof}
To save space, we write $\mathcal{G}:= \Gproj\L$ throughout the proof.

(1)
Assume that  $\rad^2\mathcal{G}=0$. 
Since $\L=\mathcal{G}(\L,\L)$, we have 
\begin{align}
\rad^2\L = \rad\L\cdot\rad\L= \rad\mathcal{G}(\L,\L)\cdot\rad\mathcal{G}(\L,\L) \subseteq  \rad^2\mathcal{G}(\L,\L)=0.  
\end{align}
If $\L$ is not CM-free, then there exists an indecomposable non-projective Gorenstein-projective $\L$-module $S$, which is simple by \cite[Lemma 2.1]{X-WChen_2012}.
Choose a projective cover $\pi:P\to S$ and a minimal left $\proj\L$-approximation $f: S \to Q$ of $S$.
Note that $f$ is the kernel of a projective cover of some $S^\prime\in\ind\sGproj\L$.
We also note that both $P$ and $Q$ are uniserial modules of length two.
Then we have $0\neq f\pi \in \rad\mathcal{G}(S, Q)\cdot\rad\mathcal{G}(P,S) \subseteq  \rad^2\mathcal{G}(P,Q)  =0$, a contradiction.
Therefore, $\L$ is CM-free.
Conversely, suppose $\mathcal{G}  = \proj\L$ and $\rad^2\L=0$.
The second equality implies that $\rad\mathcal{G}(\L^n,\L^m)\cdot\rad\mathcal{G}(\L^l,\L^n) = 0$ for all $l,m,n \geq 1$.
Then the first equality concludes that  $\rad^2\mathcal{G}=0$.

(2) 
If  $\rad\mathcal{G}=0$, then $\rad \L = \rad\mathcal{G}(\L,\L)=0$.
Conversely, if $\rad\L=0$, then all indecomposable $\L$-modules are simple, so that $\rad\mathcal{G}=0$.
\end{proof}



\subsection{Main results} 

The aim of this subsection is to characterize CM-finite G-radical-square-zero algebras and to indicate an additive generator of $\sGproj\L$ for such $\L$.
As a byproduct, we describe the singularity category $\D_{\rm sg}(\L)$ of an Iwanaga-Gorenstein G-radical-square-zero algebra $\L$.

For a $\mathbbm{k}$-linear Hom-finite Krull-Schmidt Frobenius exact category $\C$, we denote by $\Proj(\C)$ the subcategory of projective objects in $\C$.
Its stable category $\sC$, defined similarly for $\smod\L$, is a $\mathbbm{k}$-linear Hom-finite Krull-Schmidt triangulated category such that $\ind\sC=\{X\in\ind\C \mid X \text{ is non-projective}\}$.
Recall that $u_X: S_X \to X$ (resp.~$u^X: X\to S^X$) denotes a sink morphism ending (resp.~a source morphism starting) at $X\in\ind\C$.

\begin{lemma} \label{claim_17}
Let $\C$ be a $\mathbbm{k}$-linear Hom-finite Krull-Schmidt Frobenius exact category such that $|\ind\Proj(\C)| < \infty$ and $\rad^2\sC=0$.
Then the following conditions are equivalent.
\begin{enumerate}
    \item $|\ind\C| < \infty$.
    \item Both $\Mod\C$ and $\Mod(\C^\op)$ are locally finite.
    \item $\C$ has both sink morphisms and source morphisms.
\end{enumerate}
In this case, 
\begin{align}
    \textstyle \add\left(\bigoplus_{P\in\ind\Proj(\C)} S_P \oplus \Sigma^{-1}S^P\right) 
    = \sC 
    = \add\left(\bigoplus_{P\in\ind\Proj(\C)} S_P \oplus \nu S_P\right).
\end{align}
In particular,
\begin{align}
    \textstyle 
    \sC =     
    \begin{cases}
    \add\left(\bigoplus_{P\in\ind\Proj(\C)} S_P \right) & \text{if $\sC$ is semisimple,}\\[8pt]
    \add\left(\bigoplus_{P\in\ind\Proj(\C)} S_P \oplus \Sigma^{|l|}S_P\right) & \text{if $\sC$ is $l$-Calabi-Yau.}
    \end{cases}
\end{align}
In both cases,  $\sC = \thick\left(\bigoplus_{P\in\ind\Proj(\C)} S_P \right)$.
\end{lemma} 
\begin{proof}
It suffices to prove (3) $\Rightarrow$ (1), since (1) $\Rightarrow$ (2) and (2) $\Rightarrow$ (3) follow from \cref{example_1} and \cref{claim_16}, respectively.

Assume (3) holds.
Then $\C$ has Auslander-Reiten conflations by \cite[Lemma 3.2]{Iyama-Nakaoka-Palu_2024}, so that $\sC$ has Auslander-Reiten triangles.
Let  $X \xrightarrowtail{u} Y\xtwoheadrightarrow{v}Z$ be an Auslander-Reiten conflation in $\C$.
Together with \cref{claim_3}, the induced Auslander-Reiten triangle $X\xrightarrow{\underline{u}}Y\xrightarrow{\underline{v}}Z\to \Sigma X$ in $\sC$ can be used to obtain that  $Y\in\Proj(\C)$ or $Y = Y^\prime \oplus Q$ in $\C$ for some $Y^\prime\in\ind\sC$ and $Q\in\Proj(\C)$.
If $Y\in\C$ admits a direct summand $P \in\ind\Proj(\C)$, then one can deduce as in the proof of \cite[Proposition 3.4]{Hafezi-Bahlekeh_2024} that $X$ is a direct summand of $S_P$. 
Indeed, since the inflation $u: X \rightarrowtail Y$ is a source morphism, there exists an irreducible morphism $X\to P$, and hence the existence of a sink morphism $S_P\to P$ implies that $S_P$ contains $X$ as a direct summand. 
Consider now the case $Y\in\ind\sC$.
Then \cref{claim_26} implies that the induced sink morphism $\underline{v}: Y \to Z$ in $\sC$ is also a source morphism.
This shows the existence of an Auslander-Reiten conflation in $\C$ of the form 
\begin{align}\label{eq_5}
Y\xrightarrowtail{\left[\substack{v\\v^\prime}\right]}Z\oplus Q\xtwoheadrightarrow{w^\prime} W,
\end{align}
where $Q\in\Proj(\C)$ and $W\in\ind\sC$.
Then $Q\neq 0$. 
Indeed, if $Q=0$, the deflation $v: Y \twoheadrightarrow Z$ becomes an inflation and thus an isomorphism, a contradiction.
Since there exists an irreducible morphism in $\C$ from an indecomposable summand $P$ of $Q$ to $W$, the object $W$ is a direct summand of $S^P$.
Then the induced triangle $X \xrightarrow{\ \underline{u}\ }Y \xrightarrow{\ \underline{\left[\substack{v\\v^\prime}\right]}\ }Z\oplus Q \xrightarrow{\ \underline{w^\prime}\ }W$ in $\sC$ yields that $X \cong \Sigma^{-1}W$ is a direct summand of $\Sigma^{-1}S^P$.
So far, we have shown that 
\begin{align}
\textstyle 
\sC =\add\left(\bigoplus_{P\in\ind\Proj(\C)} S_P \oplus \Sigma^{-1}S^P\right).
\end{align}
Since $|\ind\Proj(\C)|<\infty$ by assumption, it follows that  $|\ind\sC|<\infty$ and hence $|\ind\C|=|\ind\sC|+|\ind\Proj(\C)|<\infty$.

To obtain another additive generator of $\sC$, we observe that if $Y \in \ind\sC$, then \eqref{eq_5} shows that $Y$ is a direct summand of $S_P$ for some indecomposable summand $P$ of $Q$.
Since $\nu X \cong Y$ in $\sC$ by \cref{claim_26}(2), where $\nu $ is a Serre functor for $\sC$, we conclude that $X \cong \nu ^{-1}Y$ is a direct summand of $\nu ^{-1}S_P$.
Therefore, we have
\begin{align} 
\textstyle 
\sC =\add\left(\bigoplus_{P\in\ind\Proj(\C)} S_P \oplus \nu ^{-1}S_P\right)=\add\left(\bigoplus_{P\in\ind\Proj(\C)} \nu S_P \oplus S_P\right).
\end{align}
These equalities show that 
\begin{align}
\textstyle 
\sC =\add\left(\bigoplus_{P\in\ind\Proj(\C)} S_P \oplus \Sigma^{|l|}S_P\right)
\end{align}
if $\sC$ is $l$-Calabi-Yau.
Moreover, if $\sC$ is semisimple, then, since $\nu X=X$  for all $X\in\ind\sC$ by \cref{claim_5}, we have $\sC = \add\left(\bigoplus_{P\in\ind\Proj(\C)} S_P \right)$.
\end{proof}

We first consider G-radical-square-zero self-injective algebras.

\begin{proposition}\label{claim_25}
Let $\L$ be a self-injective algebra.
\begin{enumerate}
\item The following conditions are equivalent.
\begin{enumerate}
\item $\L$ is G-radical-square-zero.
\item $\sGproj\L=\add\left(\rad\L\oplus\soc\L\right)$.
\item  $\L$ is a Nakayama algebra with radical cube zero.
\end{enumerate}
\item The following conditions are equivalent.
\begin{enumerate}
\item $\L$ is G-semisimple.
\item $\sGproj\L=\add(\rad\L)$.
\item  $\L$ is a Nakayama algebra with radical square zero.
\end{enumerate} 
\end{enumerate}
\end{proposition} 
\begin{proof}

For (1), it suffices to show (a) $\Rightarrow$ (b) $\Rightarrow$ (c).

(a) $\Rightarrow$ (b)
Assume (a) holds.
Since $\L$ is self-injective, for each $P\in\ind\proj\L$, the canonical inclusion $\rad P \hookrightarrow P$ and the canonical surjection $P\to P/\soc P$ are a sink morphism and a source morphism in $\Gproj\L$, respectively.
Thus by \cref{claim_17},  
\begin{align}
\textstyle 
    \sGproj\L =\add\left(\bigoplus_{P\in\ind\proj\L} S_P \oplus \Sigma^{-1}S^P\right) = \add \left(\rad\L\oplus\soc\L\right).
\end{align}

(b) $\Rightarrow$ (c)
Assume  $\Gproj\L = \add \left(\rad\L\oplus\soc\L\oplus\L\right).$
We claim that every indecomposable projective $\L$-module $P$ is a uniserial module of length at most three.
This immediately implies (c) (since $\L$ is self-injective).
It is clear that if $\rad P$ is simple, then $P$ is a uniserial module with $l(P)=2$.  
Suppose that $\rad P$ is non-simple, and choose a maximal submodule $M$ of $\rad P$.
Since $\rad P$ has simple socle, $M$ is indecomposable.
Moreover, since $M\in\ind\sGproj\L$, we have $M=\rad Q$ or $M=\soc Q$ for some $Q\in\ind\proj\L$ by assumption.
If $M=\rad Q$, then the fact that $P$ and $Q$ have the same socle implies that $P=Q$, a contradiction.
Thus $M=\soc Q$ is simple, and hence $P$ is uniserial with $l(P)=3$.

For (2), it suffices to show (b) $\Rightarrow$ (c), since (a) $\Rightarrow$ (b) follows form \cite[Proposition 3.13]{Hafezi-Bahlekeh_2024}.
Assume (b) holds.
By (1), $\L$ is a self-injective Nakayama algebra with $\rad^3\L=0$.
Thus $\Gproj\L=\add(\L\oplus\rad\L\oplus\rad^2\L)$.
However, since $\Gproj\L=\add(\L\oplus\rad\L)$ by assumption, we must have $\rad^2\L=0$.
\end{proof}

We denote by $X_M=X_{M_\L}$  a minimal right $\Gproj\L$-approximation of $M\in\mod\L$ (if it exists).
The following is the main result of this section.

\begin{theorem}\label{claim_19}
The following conditions are equivalent for a G-radical-square-zero algebra $\L$.
\begin{enumerate}
    \item $\L$ is CM-finite.
    \item $\Gproj\L$ is contravariantly finite in $\mod\L$.
\end{enumerate}
In this case, 
\begin{align}
    \textstyle 
    \add\left(X_{\rad\L}\oplus \Omega_{\L}(X_{\rad\L_{\L^\op}}^*)\right) 
    = \sGproj\L 
    = \add\left(X_{\rad\L} \oplus \nu X_{\rad\L}\right).
\end{align}
In particular,
\begin{align}
    \textstyle 
    \sGproj\L  =     
    \begin{cases}
    \add(X_{\rad\L}) & \text{if $\L$ is G-semisimple,}\\[5pt]
    \add\left(X_{\rad\L}\oplus \Omega_\L^{|l|}(X_{\rad\L})\right) & \text{if $\sGproj\L $ is $l$-Calabi-Yau.}
    \end{cases}
\end{align}
In both cases, $\sGproj\L = \thick(X_{\rad\L})$.
\end{theorem} 
\begin{proof}
(1) $\Rightarrow$ (2) 
Follows from \cite[Proposition 4.2]{Auslander-Smalo_1980}.

(2) $\Rightarrow$ (1)
Assume (2) holds.
Since $\Gproj\L$ is functorially finite in $\mod\L$ by \cite[Corollary 2.6]{Krause-Solberg_2003}, it follows from \cite[Theorem 2.4]{Auslander-Smalo_1981} that $\Gproj\L$ has almost split sequences in the sense of \cite[p.429]{Auslander-Smalo_1981}.
Since, from the definition, $\Gproj\L$ has sink morphisms and source morphisms,  \cref{claim_17} implies that $|\ind\Gproj\L| <\infty$.
For the remaining, assume that $\L$ is CM-finite.
Then by \cite[Lemma 3.2]{Hafezi-Bahlekeh_2024}, we have $\bigoplus_{P\in\ind\proj\L} S_P = X_{\rad\L}$ and $\bigoplus_{P\in\ind\proj\L} S^P = X_{\rad\L_{\L^\op}}^*$.
Further, if $\sGproj\L$ is $l$-Calabi-Yau, then one of $\nu=\Sigma^{l}$ and $\nu^{-1}=\Sigma^{-l}$ coincides with $\Omega_\L^{|n|}$, where we use the fact that $\Sigma^{-1}=\Omega_\L$.
Together with these observations,  \cref{claim_17} yields the desired descriptions of $\sGproj\L$.
\end{proof}

\begin{example}\label{example_8}
Consider the upper triangular matrix algebra $T_2(\L) = \begin{bmatrix}
        \L & \L \\
        0 & \L \\
\end{bmatrix}$ 
of a finite dimensional CM-finite G-semisimple $K$-algebra $\L$. 
Then \cite[Proposition 3.7, Lemma 5.5 and Theorem 5.8]{Hafezi-Bahlekeh_2024} imply that $T_2(\L)$ is a CM-finite non-G-semisimple G-radical-square-zero algebra such that every path-closed subset of $\ind\sGproj T_2(\L)$ has $3n$ elements for some $n\geq 1$. 
In particular, $\sGproj T_2(\L)$ cannot be Calabi-Yau because of \cref{claim_42,claim_36}.
\end{example}

\begin{remark}
In \cref{example_3}, we give an example of a CM-finite non-G-semisimple G-radical-square-zero algebra $\L$ such that $\sGproj\L$ is Calabi-Yau.
\end{remark}

We say that $M\in\mod\L$ is {\it eventually periodic} if $\Omega_\L^{n+p}(M) \cong \Omega_\L^{n}(M)$ in $\mod\L$ for some $n\geq 0$ and $p \geq 1$.
If $n=0$, we call $M$ {\it periodic}.
In this case, the least $p \geq 1$ with the property that $\Omega_\L^{p}(M) \cong M$ in $\mod\L$ is called the {\it period} of $M$.
If $M$ is eventually periodic, then all periodic higher syzygies $\Omega_\L^n(M)$ have the same period $p$.
In this case, we call $M$ an {\it eventually periodic module of period $p$}.

\begin{lemma}\label{claim_20_lemma}
    Let $\L$ be a $d$-Iwanaga-Gorenstein algebra.
    Then any eventually periodic $\L$-module $M$ of period $p$ admits a right $\Gproj\L$-approximation of the form $\Omega_\L^n(M) \oplus P \to M$ for some $d \leq n < d+p$ and $P\in \proj\L$.
\end{lemma} 
\begin{proof}
Let $M$ be an eventually periodic $\L$-module of period $p$.
Since the $d$-th syzygy $\Omega_\L^d(M) \in \Gproj\L$  satisfies $\Omega_\L^d(M) \cong \Omega_\L^{d+p}(M)$ in $\sGproj\L$ (cf.~\cite[Theorem 3.9]{Usui_2024}), there exists an exact sequence
\begin{align}
    0 \to \Omega_\L^{d}(M) \to P_{d+p-1} \to \cdots \to P_{d} \to \Omega_\L^{d}(M) \to 0
\end{align} 
in $\Gproj\L$ with $P_i\in\proj\L$.
Then splicing copies of this exact sequence yields an acyclic complex $T_\bullet$ of projective modules.
Let $f: X_M \to M$ be a minimal right $\Gproj\L$-approximation of $M$.
Since $T_\bullet$ can be used to compute $\Sigma^i(\Omega_\L^{j}(M))$ for all $i \in \Z$ and $j \geq d$, there exists an integer $n$ with $d \leq n< d+p$ such that $X_M \cong \Sigma^d(\Omega_\L^{d}(M)) \cong \Omega_\L^{n}(M)$ in $\sGproj\L$, where the first isomorphism follows from \cite[Lemma 3.6]{Auslander-Reiten_1991}.
Choose $P, Q \in \proj\L$ such that $X_M \oplus Q \cong \Omega_\L^{n}(M) \oplus P$ in $\Gproj\L$.
Then the composition
\begin{align}
    \Omega_\L^{n}(M) \oplus P \xrightarrow{\sim} X_M \oplus Q \xrightarrow{[f\ 0]} M
\end{align}
is the desired right $\Gproj\L$-approximation of $M$.
\end{proof}

Over a CM-finite Iwanaga-Gorenstein algebra $\L$, every $\L$-module is eventually periodic (use \cite[Lemma 2.1.13 and Theorem 2.3.7]{X-WChen_2017} for example).
We then set
\begin{align}\label{eq_13}
    n_{\L} := \min\{ n \mid d-1 \leq n, \,  X_{\rad\L} \cong \Omega_\L^{n}(\rad\L) \text{ in } \sGproj\L \},
\end{align}
where $d=\injdim\L$. 
We now specialize \cref{claim_19} to Iwanaga-Gorenstein algebras.

\begin{proposition}\label{claim_20}
Let $\L$ be a $d$-Iwanaga-Gorenstein G-radical-square-zero algebra.
Then $\L$ is CM-finite, and 
\begin{align}
    \textstyle 
    \sGproj\L 
    &=\add\left(\Omega_\L^{d}(\L/\rad\L) \oplus \Sigma^{n_{\L}+n_{\L^\op}-2d+1}\hspace{1.5pt}\Omega_{\L^\op}^{d}\!\left(\L/\rad\L\right)^*\right)\\[1pt]
    &= \add\left(\Omega_\L^{d}(\L/\rad\L) \oplus \nu \Omega_\L^{d}(\L/\rad\L)\right),
\end{align}
where $n_{\L}$ and $n_{\L^\op}$ are defined in \eqref{eq_13}.
In particular,
\begin{align}
    \textstyle 
    \sGproj\L =     
    \begin{cases}
    \add(\Omega_\L^{d}(\L/\rad\L)) & \text{if $\L$ is G-semisimple,}\\[1pt]
    \add\left(\Omega_\L^{d}(\L/\rad\L) \oplus \Omega_\L^{d+|l|}(\L/\rad\L)\right) & \text{if $\sGproj\L$ is $l$-Calabi-Yau.}
    \end{cases}
\end{align}
In both cases, $\sGproj\L = \thick(\Omega_\L^{d}(\L/\rad\L))$.
\end{proposition} 
\begin{proof}
The CM-finiteness of $\L$ follows from \cref{claim_19} and the contravariant finiteness of $\Gproj\L$ (see \cite[Corollary 2.3.6]{X-WChen_2017} for example).
Moreover, \cref{claim_20_lemma} implies that
\begin{align}\label{eq_8}
X_{\rad\L} \in \add(\Omega_\L^{n_\L}(\rad\L)) \subseteq\sGproj\L\quad\text{and}\quad    X_{\rad\L_{\L^\op}} \in \add(\Omega_{\L^\op}^{n_{\L^\op}}(\rad\L))\subseteq\sGproj\L^\op.
\end{align}
Note that 
\begin{align}
    \Omega_\L^{n_\L}(\rad\L) 
    = \Omega_\L^{n_\L-(d-1)}\!\left(\Omega_\L^{d-1}\!\left(\rad\L\right)\right)
    = \Sigma^{-n_\L+d-1}\Omega_\L^{d}\!\left(\L/\rad\L\right)
\end{align}
and similarly $\Omega_{\L^\op}^{n_{\L^\op}}(\rad\L) = \Sigma^{-n_{\L^\op}+d-1}\Omega_{\L^\op}^{d}\!\left(\L/\rad\L\right)$.
Thus we have 
\begin{align}
    \sGproj\L 
    &= \add\left(X_{\rad\L} \oplus \Omega_{\L}(X_{\rad\L_{\L^\op}}^*) \right) \hspace{36mm} (\text{by \cref{claim_19}})\\
    &= \add\left(\Omega_\L^{n_\L}(\rad\L) \oplus \Omega_{\L}\!\left(\Omega_{\L^\op}^{n_{\L^\op}}(\rad\L)^*\right) \right)  \hspace{20.5mm} (\text{by \eqref{eq_8}})\\
    &= \add\left(\Sigma^{-n_\L+d-1}\hspace{1pt}\Omega_\L^{d}\left(\L/\rad\L\right) \oplus \Sigma^{-1}\left(\Sigma^{-n_{\L^\op}+d-1}\Omega_{\L^\op}^{d}\left(\L/\rad\L\right)\right)^*\right)\\
    &= \add\left(\Sigma^{-n_\L+d-1}\hspace{1pt}\Omega_\L^{d}\left(\L/\rad\L\right) \oplus \Sigma^{n_{\L^\op}-d}\hspace{1.5pt}\Omega_{\L^\op}^{d}\left(\L/\rad\L\right)^* \right)\\
    &= \add\left(\Omega_\L^{d}\left(\L/\rad\L\right) \oplus \Sigma^{n_\L+n_{\L^\op}-2d+1}\hspace{1.5pt}\Omega_{\L^\op}^{d}\left(\L/\rad\L\right)^* \right).
\end{align}
The other descriptions of $\sGproj\L$  also follow from \cref{claim_19} and \eqref{eq_8}.
\end{proof}

Thanks to Buchweitz \cite{Buchweitz_1986}, we obtain the following consequence of \cref{claim_20}.

\begin{corollary}\label{claim_41}
Let $\L$ be an Iwanaga-Gorenstein  G-radical-square-zero algebra.
Then 
\begin{align}
    \textstyle 
    \add\left(\L/\rad\L \oplus \Sigma^{1+n_{\L}+n_{\L^\op}}(\L/\rad\L)^*\right)=\D_{\rm sg}(\L)
    = \add\left(\L/\rad\L \oplus \nu\L/\rad\L\right),
\end{align}
where $n_{\L}$ and $n_{\L^\op}$ are defined in \eqref{eq_13}.
In particular,
\begin{align}
    \textstyle 
    \D_{\rm sg}(\L)=     
    \begin{cases}
    \add(\L/\rad\L) & \text{if $\D_{\rm sg}(\L)$ is semisimple,} \\[1pt]
    \add\left(\L/\rad\L \oplus \Sigma^l\L/\rad\L\right) & \text{if $\D_{\rm sg}(\L)$ is $l$-Calabi-Yau.}
    \end{cases}
\end{align}
In both cases,  $\D_{\rm sg}(\L) = \thick(\L/\rad\L)$.
\end{corollary} 
\begin{proof}
Observe that $\Omega_\L^i(X) \cong \Sigma^{-i}X$ in $\D_{\rm sg}(\L)$ for $X\in\mod\L$ and $i\geq 0$ (\cite[Lemma 2.2.1]{Buchweitz_1986}).
Since we have  $\Sigma^{n_\L+n_{\L^\op}-2d+1}\hspace{1.5pt}\Omega_{\L^\op}^{d}\!\left(\L/\rad\L\right)^*  \cong \Sigma^{1-d+n_\L+n_{\L^\op}}\!\left(\L/\rad\L\right)^* $ in $\D_{\rm sg}(\L)$, the statements follow from \cref{claim_20} and \cite[Theorem 4.4.1]{Buchweitz_1986}.
\end{proof}

Two artin algebras $\L$ and $\Gamma$ are called {\it singularly equivalent} if there is a triangle equivalence between their singularity categories $\D_{\rm sg}(\L)$ and $\D_{\rm sg}(\Gamma)$.
It is known that the property of being Iwanaga-Gorenstein is not preserved under singular equivalences; see \cite[Example 4.3]{X-WChen_2009} for example.

\begin{example}\label{example_6}
Let $\L$ and $\G$ be singularly equivalent algebras.
If $\L$ is an Iwanaga-Gorenstein G-radical-square-zero algebra, then $\G$ is G-radical-square-zero (since $\sGproj\G$ is a thick subcategory of the radical-square-zero $\mathbbm{k}$-linear Hom-finite Krull-Schmidt triangulated category $\D_{\rm sg}(\G)$).
Moreover, since 
\begin{align}
    |\ind\sGproj\G| \leq |\ind\D_{\rm sg}(\G)| = |\ind\D_{\rm sg}(\L)| < \infty,
\end{align}
\cref{claim_19} gives an additive generator of $\sGproj\G$. 
\end{example}

We now focus on G-semisimple algebras.
Following \cite[Definition 3.5]{Hafezi-Bahlekeh_2024}, we say that $X\in\ind\sGproj\L$ is {\it G-semisimple} if the short exact sequence $0 \to \Omega_\L(X)  \to P \xrightarrow{\pi} X \to 0$ with $\pi$ a projective cover of $X$ is an Auslander-Reiten conflation in $\Gproj\L$.

\begin{lemma}\label{claim_23}
The following conditions are equivalent for an artin algebra $\L$.
\begin{enumerate}
\item $\L$ is G-semisimple.
\item Every indecomposable non-projective object of $\Gproj\L$ is G-semisimple.
\end{enumerate}
\end{lemma} 
\begin{proof}
It is a consequence of \cref{claim_6}, since $X\in\ind\sGproj\L$ is G-semisimple if and only if the middle term $S_X$ of an Auslander-Reiten conflation in $\Gproj\L$ ending at $X$ is projective.
\end{proof}

\begin{remark}
The proposition recovers \cite[Lemma 3.6]{Hafezi-Bahlekeh_2024}.
\end{remark}

It was proved in \cite[Proposition 3.13]{Hafezi-Bahlekeh_2024} that an artin algebra $\L$ is $1$-Iwanaga-Gorenstein G-semisimple if and only if $\Gproj\L = \add (\L \oplus \rad\L)$. 
We generalize this result as follows.

\begin{proposition}\label{claim_35}
Let $\L$ be an artin algebra such that every irreducible morphism in $\Gproj\L$ between indecomposable non-projective objects is either a monomorphism or an epimorphism.
Then the following conditions are equivalent.
\begin{enumerate}
\item $\L$ is a CM-finite G-semisimple algebra.
\item $\Gproj\L = \add(\L \oplus X_{\rad\L})$.
\end{enumerate}
\end{proposition}
\begin{proof}
(1) $\Rightarrow$ (2) 
Follows from \cref{claim_19}.

(2) $\Rightarrow$ (1)
Assume (2) holds.
Clearly, $\L$ is CM-finite.
Thus $\Gproj\L$ has Auslander-Reiten conflations (cf.~the proof of \cref{claim_19}).
The argument below parallels that in the proof of \cite[Proposition 3.13]{Hafezi-Bahlekeh_2024}.
We claim that $\Gproj\L$ has no irreducible morphism between indecomposable non-projective objects.
Assume to the contrary that there is such an irreducible morphism $f: G \to G^\prime$ in $\Gproj\L$.
Then $f$ is either a monomorphism or an epimorphism by assumption.
Assume first that $f$ is an epimorphism.
Then the middle term $S_{G^\prime} = S^{\tau G^\prime}$ of an Auslander-Reiten conflation ending at $G^\prime$ must contains an indecomposable projective object $P\in\ind\proj\L$, since $\tau G^\prime$ is a direct summand of $X_{\rad\L} = \bigoplus_{P\in\ind\proj\L}S_P$.
Thus we obtain an irreducible morphism $g : P \to G^\prime$ in $\Gproj\L$.
Since $\Gproj\L$ is closed under kernels of epimorphisms,  any epimorphism in $\Gproj\L$ is a deflation. 
Hence $g : P \to G^\prime$  factors through the irreducible deflation $f: G \to G^\prime$, a contradiction.
Next, consider the case where $f$ is a monomorphism.
Since $\L$ is CM-finite, there is an integer $n>0$ such that $\tau^n G^\prime \cong G^\prime$.
For each $1\leq i \leq n$, take an Auslander-Reiten conflation $0\to \tau^i G^\prime \to X_i \to \tau^{i-1}G^\prime \to 0$ in $\Gproj\L$.
Note that $G$ is a direct summand of $X_1$.
Then there is a chain 
\begin{align}
    G^\prime = \tau^n G^\prime \to M_n \to \cdots \to \tau^2 G^\prime \to M_2 \to \tau G^\prime \to G \xrightarrow{f} G^\prime
\end{align}
of irreducible morphisms in $\Gproj\L$ between indecomposable non-projective objects, where each $M_i$ is a direct summand of $X_i$.
Since irreducible morphisms in $\Gproj\L$ are all monomorphisms, we conclude that $l(G) < l(M_n) < \cdots < l(G^\prime) < l(G)$, a contradiction.
Now, the established claim implies that every Auslander-Reiten conflation in $\Gproj\L$ has projective middle term.
Therefore, $\L$ is G-semisimple by \cref{claim_23}.
\end{proof}

\begin{remark}
Monomial algebras (\cref{claim_38}) and $1$-Iwanaga-Gorenstein algebras satisfy the assumption on $\L$ in the proposition.    
In particular, it recovers \cite[Proposition 3.13]{Hafezi-Bahlekeh_2024}.
\end{remark}




\subsection{The case of monomial algebras} 
\label{The_case_of_monomial_algebras}

The aim of this subsection is to characterize  G-radical-square-zero monomial algebras.
We first recall some necessary notions from \cite{Chen-Shen-Zhou_2018, Honma-Usui_2025}.
See \cite{ASS06,ARS_1995} for the basics of bound quiver algebras.

Let $Q$ be a finite quiver, that is, a quiver such that the vertex set $Q_0$ and the arrow set $Q_1$ are both finite.
For $n\geq1$, a {\it path of length $n$} in $Q$ is a sequence $p=a_1\cdots a_n$ of arrows $a_i \in Q_1$ such that the target $t(a_i)$ of $a_i$ equals the source $s(a_{i+1})$ of $a_{i+1}$ for $1\leq i< n$.
Define $s(p):=s(a_1)$ and $t(p):=t(a_n)$.
The length of $p$ is denoted by $l(p)$.
For each $i\in Q_0$, the corresponding trivial path is denoted by $e_i$. 
Note that $s(e_i)=i=t(e_i)$. 
We denote by $\B$ the set of paths
and by $\B_{\geq i}$ the set of paths of length at least $i$. 
Note that $\B_{+} := \B_{\geq 1}$ is the set of non-trivial paths.
For $p, q \in \B$, we say that $q$ is a {\it subpath} of $p$ if $p = p^{\prime} q p^{\prime\prime}$ for some $p^{\prime}, p^{\prime\prime} \in \B$.
If $p^{\prime}$ (resp.~$p^{\prime\prime}$) is trivial, we call $q$ a {\it left \text{\rm (resp.}~right\text{\rm )} divisor} of $p$. 
For a subset $S\subseteq\B$, we say that $p\in S$ is {\it minimal in $S$} if any subpath of $p$ lying in $S$ equals $p$.
We call a non-trivial path $c$ with $s(c) = t(c)$ a {\it cycle}.
A cycle $c$ is called {\it repetition-free} if any cycle $c^\prime$ such that $c= (c^\prime)^m$ for some $m \geq 1$ equals $c$.

In the rest of this section, let $\L = kQ/I$ be a monomial algebra, that is, a bound quiver algebra such that the admissible ideal $I$ is generated by paths. 
We denote by $\F$ the set of minimal paths in $\B\cap I$.
This set is finite, since $I$ is finitely generated.
For $p \in \B$, we simply write the canonical image $p+I$ in $\L$ as $p$ when there is no danger of confusion.
We say that $p \in \B$ is {\it non-zero in $\L$} if $p\neq0$ in $\L$ (i.e.~$p\not\in I$).

Recall from \cite[Definitions 3.3 and 3.7]{Chen-Shen-Zhou_2018} that a pair $(p, q)$ of non-zero non-trivial paths is {\it perfect} if the following conditions hold.
\begin{enumerate}
\item[(i)] $t(p) = s(q)$ and $pq\in I$.
\item[(ii)] If $pq^\prime\in I$ for a non-zero path $q^\prime$ with $t(p) = s(q^\prime)$, then $q$ is a left divisor of $q^\prime$.
\item[(iii)] If $p^\prime q\in I$ for a non-zero path $p^\prime$ with $t(p^\prime) = s(q)$, then $p$ is a right divisor of $p^\prime$.
\end{enumerate}
A sequence $(p_1, \ldots, p_n, p_{n+1}=p_1)$ of non-zero non-trivial paths is called a {\it perfect path sequence} if $(p_i, p_{i+1})$ is a perfect pair for all $ 1 \leq i \leq n$. 
In this case,  each $p_i$ is called a {\it perfect path}. 
Let $\mathbb{P}_{\L}$ be the set of perfect paths in $\L$.
This set is finite, since $\L$ has only finitely many non-zero non-trivial paths. 
By \cite[Theorem 4.1]{Chen-Shen-Zhou_2018}, there is a bijection between $\mathbb{P}_{\L}$ and $\ind \sGproj\L$, where we assign to each $p\in\P_\L$ the cyclic $\L$-module $p\L$.
Therefore, $\L$ is CM-finite, and $\L$ is CM-free if and only if $\mathbb{P}_\L$ is empty.

For each $p_1\in\mathbb{P}_\L$, choose a perfect path sequence $(p_{1}, \ldots, p_{n}, p_{n+1} = p_{1})$.
The repetition-free subcycle $c$ of the cycle $p_1 \cdots p_n$ satisfying $p_1 \cdots p_n = c^l$ for some $l>0$ is called the {\it underlying cycle associated with $p_1$} and denoted by $c_{p_1}$. 
Note that $c_{p_1}$ is independent of the choice of a perfect path sequence.
We denote by $\C(\L)$ the set of  equivalence classes (with respect to cyclic permutation) of underlying cycles.
We call $p\in\P_\L$ {\it elementary} (resp.~{\it co-elementary}) if there is no $q\in\P_\L$ such that $q=pr$ (resp.~$p=qr$) for some $r\in\B_+$.
Denote by $\E_\L$ (resp.~$\Eco_\L$) the set of elementary (resp.~co-elementary) paths.
By \cite[Theorem 2.11]{Honma-Usui_2025} (resp.~\cite[Proposition-Definition 2.12]{Honma-Usui_2025}), every $p\in\P_\L$ (resp.~$c\in\C(\L)$) can be written uniquely as a product $p = r_1\cdots r_n$ (resp.~$c = r_1\cdots r_n$) of co-elementary paths $r_1, \ldots, r_n \in \Eco_\L$.
Denote $|p|:=n$ (resp.~$|c|:=n$).
Note that $c\in\C(\L)$ is perfect if and only if $c\not\in I$.
For each $c=r_1 \cdots r_n \in \C(\L)$ with $r_i\in\Eco_\L$, define 
\begin{align}
    \P_\L(c) :=\{p\in\P_\L\mid p= r_1 p^\prime \text{ for some $p^\prime\in\B$}\}.
\end{align}
Moreover, define
\begin{align}
    \Fp :=\{ f \in \F \mid f=pq \text{ for some $p, q \in \P_\L$} \}.
\end{align}
We call elements of $\Fp$ {\it perfect relations}.
For each $f=r_1\cdots r_n \in \Fp$ with $r_i\in\Eco_\L$, denote $|f|:= n$. 
By \cite[Lemma 3.4]{Honma-Usui_2025}, if  $f, f^\prime \in\Fp$ are perfect relations lying on the same underlying cycle $c \in \C(\L)$, then $|f| = |f^\prime|$.
We denote by $r_c$ the common value and call it the {\it length of perfect relations on $c$}. 
For $c\in\C(\L)$ and $p\in\E_\L$ with $c_p = c$ in $\C(\L)$, we have
\begin{align}\label{monomial_eq_1}
    |p| = |\P_\L(c)| = r_c -1.
\end{align}

Let $p=rq \in \P_\L$ with $r \in \Eco_\L$ and $q\in\B$ be an arbitrary perfect path.
Note that $q\in\P_\L$ whenever $q\in\B_+$ (\cite[Proposition 2.5]{Honma-Usui_2025}). 
Choose a co-elementary subpath $s$ of $c_p$ such that $qs\in\P_\L$.
Then there is a short exact sequence 
\begin{align}\label{eq_9}
    0 \to qs\L \xrightarrow{\,\csubstack{\iota^\prime}{-\pi^\prime}\,} q\L \oplus ps\L \xrightarrow{\,\rsubstack{\pi}{\iota}\,} p\L \to 0
\end{align}
in $\Gproj\L$, where $\iota^\prime: qs\L \hookrightarrow q\L$ and $\iota: ps\L \hookrightarrow p\L$ are the canonical inclusions, and $\pi^\prime: qs\L \to rqs\L=ps\L$ and $\pi: q\L \to rq\L =p\L$
are the left multiplications by $r$.
An argument as in the proof of \cite[Proposition 3.10]{Honma-Usui_2025} proves the following proposition.
The proof is left to the reader.

\begin{proposition}
Let $\L$ be a monomial algebra.
Then the short exact sequence \eqref{eq_9} is an Auslander-Reiten conflation in $\Gproj\L$. 
\end{proposition}

We then obtain the following consequence.

\begin{corollary}\label{claim_38}
Let $\L$ be a monomial algebra.
Then any irreducible morphism in $\Gproj\L$ between indecomposable non-projective objects is either a monomorphism or an epimorphism.
\end{corollary}

It was proved in \cite[Theorem 4.3]{Honma-Usui_2025} that there is a triangle equivalence
\begin{align}
 \sGproj\L\, \cong \prod_{c \in \C(\L)} \D^{\rm b}\!\left(\mod k A_{r_c -1}\right)\!\big/\tau^{|c|},   
\end{align}
where 
$\tau$ denotes the Auslander-Reiten translation for $\D^{\rm b}\!\left(\mod k A_{r_c -1}\right)$. 
Therefore, we have an isomorphism   
\begin{align}\label{monomial_eq_2}
 \G(\sGproj\L)\, \cong \prod_{c \in \C(\L)} \Z\!A_{r_c -1}\big/\tau^{|c|} 
\end{align}
of stable valued translation quivers.
The following proposition characterizes G-radical-square-zero monomial algebras in terms of perfect paths.

\begin{proposition}\label{claim_27}
The following conditions are equivalent for a monomial algebra $\L$.
\begin{enumerate}
\item $\L$ is G-radical-square-zero.
\item Every $p\in\P_\L$ satisfies $1 \leq |p| \leq 2$. 
\item Every $p\in\P_\L$ is co-elementary or elementary.
\end{enumerate}
\end{proposition}
\begin{proof} 
(1) $\Leftrightarrow$ (2)
Follows from \eqref{monomial_eq_1}, \eqref{monomial_eq_2} and \cref{claim_3}.

(2) $\Leftrightarrow$ (3)
The fact that there is a perfect path that is neither elementary nor co-elementary if and only if there is an elementary path $p\in\E_\L$ such that $|p|\geq 3$  yields the equivalence (2) $\Leftrightarrow$ (3), where, to obtain the implication $(\Leftarrow)$, we use the fact that every perfect path is a subpath of some elementary path.
\end{proof}

Following \cite[p.1130]{Chen-Shen-Zhou_2018}, we say that {\it there is no overlap in $\L$} if there is no pair of  perfect paths $p, q$ satisfying one of the following mutually exclusive conditions:
\begin{enumerate}
    \item [(O1)] $p=q$, and $p=p^\prime x$ and $q=x q^\prime$ for some $x, p^\prime, q^\prime \in \B_+$ such that $p^\prime x q^\prime \not\in I$.
    \item [(O2)] $p\not=q$, and $p=p^\prime x$ and $q=x q^\prime$ for some $x \in \B_+$ and some $p^\prime, q^\prime \in \B$ such that  $p^\prime x q^\prime\not\in I$.
\end{enumerate}
Chen, Shen and Zhou \cite[Proposition 5.9]{Chen-Shen-Zhou_2018} and the authors \cite[Proposition 4.4]{Honma-Usui_2025} already obtained characterizations of G-semisimple monomial algebras.
For latter use, we restate these results in the following proposition.

\begin{proposition}\label{claim_28}
The following conditions are equivalent for a monomial algebra $\L$.
\begin{enumerate}
\item $\L$ is G-semisimple.
\item There is no overlap in $\L$.
\item Every $p\in\P_\L$ is co-elementary.
\item Every $p\in\P_\L$ is elementary.
\end{enumerate}
In this case,  $\sGproj\L$ is $0$-Calabi-Yau.
\end{proposition}
\begin{proof}
(1) $\Leftrightarrow$ (2)
Follows from \cite[Proposition 5.9]{Chen-Shen-Zhou_2018}.

(2) $\Leftrightarrow$ (3)
Follows from \cite[Proposition 4.4]{Honma-Usui_2025}.

(3) $\Leftrightarrow$ (4)
Follows from the observation that both (3) and (4) are equivalent to the condition that $rs=0$ in $\L$ for all $r, s\in\Eco_\L$.

The last statement follows from \cref{claim_12}.
\end{proof}

In \cite[Lemma 3.3]{Elsener_2021}, Elsener characterized when  $\sGproj\L$ with $\L$ a $1$-Iwanaga-Gorenstein monomial algebra  is $3$-Calabi-Yau in terms of underlying cycles and the lengths of perfect relations.
The following lemma extends her result to arbitrary monomial algebras $\L$ and arbitrary odd integers $2l+1$.
Although the proof parallels that of \cref{claim_42}, we include it for completeness.

\begin{lemma}\label{claim_34}
The following conditions are equivalent for a monomial algebra $\L$ and an integer $l$.
\begin{enumerate}
\item $\sGproj\L$ is $(2l+1)$-Calabi-Yau.
\item For each $c\in\C(\L)$, we have  $l r_c \equiv -1 \ (\operatorname{mod} |c|)$.
\end{enumerate}
\end{lemma}
\begin{proof}
For $c\in\C(\L)$, we denote by $\Sigma_c$ (resp.~$\tau_c$) the shift functor (resp.~the Auslander-Reiten translation) for $\D^{\rm b}\!\left(\mod k A_{r_c -1}\right)$.
It follows that $\Sigma_c^2 \cong \tau_c^{-r_c}$ (cf.~\cite[8.3]{Keller_2005}).
It also follows that $\Sigma_c$ and $\tau_c$ induce the shift functor $\underline{\Sigma_c}$ and  the Auslander-Reiten translation $\underline{\tau_c}$ for $\T_c := \D^{\rm b}\!\left(\mod k A_{r_c -1}\right)/\tau^{|c|}$, respectively.

(1) $\Rightarrow$ (2)
Assume (1) holds.
Since each $\T_c$ is also $(2l+1)$-Calabi-Yau, we have 
\begin{align}
\underline{\Sigma_c}\,\underline{\tau_c} \cong \underline{\Sigma_c}^{2l+1}
\quad\text{and hence}\quad
\underline{\tau_c}^{-lr_c-1}\cong \underline{\Sigma_c}^{2l}\,\underline{\tau_c}^{-1} \cong 1_{\T_c}
\end{align}
as auto-equivalences of $\T_c$.
Since the group defining $\T_c$ is the cyclic group generated by $\tau^{|c|}$, for each $X\in\D^{\rm b}\!\left(\mod k A_{r_c -1}\right)$, there is an integer $b=b_{c, X}\in\Z$ such that  $\tau^{lr_c+1}X \cong \tau^{|c| b}X$ in $\D^{\rm b}\!\left(\mod k A_{r_c -1}\right)$.
This implies that $lr_c+1 = |c| b$, as desired.

(2) $\Rightarrow$ (1)
Assume (2) holds.
For each $c\in\C(\L)$, since $l r_c =b|c| -1$ for some $b\in \Z$, we have isomorphisms
\begin{align}
\underline{\Sigma_c}^{2l}\,\underline{\tau_c}^{-1} \cong \underline{\tau_c}^{-lr_c-1} \cong \underline{\tau_c}^{b_c|c|} \cong (\underline{\tau_c}^{|c|})^{-b_c} \cong 1_{\T_c}
\end{align}
as auto-equivalences of $\T_c$.
This shows that the Serre functor $\underline{\Sigma_c}\,\underline{\tau_c}$ for $\T_c$ is isomorphic to $\underline{\Sigma_c}^{2l+1}$.
Therefore, $\sGproj\L \cong \prod_{c \in \C(\L)} \T_c$ is $(2l+1)$-Calabi-Yau.
\end{proof}

\begin{example}\label{example_4}
Consider the Nakayama algebra $\L=kQ/I$ given by the following quiver with relations:
\begin{align}
\raisebox{8mm}{
    \xymatrix@R=9mm@C=9mm{
    1 \ar[r]^-{a_1} & 2 \ar[d]^-{a_2}\\
    4 \ar[u]^-{a_4} & 3 \ar[l]^-{a_3}\\
    }}
    \qquad a_1 a_2 a_3 a_4 = a_4 a_1 a_2 a_3 =0. 
\end{align}
Then $\P_\L=\{ a_1 a_2 a_3, a_4 \} = \Eco_\L = \E_\L$, 
and $\L$ is $2$-Iwanaga-Gorenstein.
Thus  \cref{claim_20,claim_41,claim_28} imply that $\L$ is G-semisimple, and 
\begin{align}
\sGproj\L=\add(\Omega_\L^{2}(\L/\rad\L))=\thick(\Omega_\L^{2}(\L/\rad\L))
\end{align}
and
\begin{align}
 \D_{\rm sg}(\L)=\add(\L/\rad\L)=\thick(\L/\rad\L).
\end{align}
Note that 
\begin{align}
\Omega_\L^{2}(\L/\rad\L) = a_2 a_3 a_4\L \oplus a_1 a_2 a_3\L
\end{align}
where $a_2 a_3 a_4\L \cong a_4\L$ by \cite[Proposition 4.6]{Chen-Shen-Zhou_2018}. 
\end{example}


\begin{example}\label{example_3}
Consider the Nakayama algebra $\L=kQ/I$ given by the following quiver with relations:
\begin{align}
    \xymatrix@C=10mm{
    1 \ar@<0.6ex>[r]^-{a_1} & 2 \ar@<0.6ex>[l]^-{a_2} 
    }
    \qquad a_1 a_2 a_1 a_2 a_1 a_2 =0. 
\end{align}
Then we have 
\begin{align}
    & \P_\L=\{ a_1 a_2, \  a_1 a_2 a_1 a_2\}, \quad
    \Eco_\L=\{ a_1 a_2\},\quad
    \E_\L=\{ a_1 a_2 a_1 a_2\},\\[1mm]
    &\Fp=\{a_1 a_2 a_1 a_2 a_1 a_2\},\quad
    \C(\L)=\{ c = a_1 a_2\}
\end{align}
and $\L$ is $2$-Iwanaga-Gorenstein.
Since $|c| =1$ and $r_{c}=3$, the integers $l=0$ and $b_{c}=1$ satisfy $l r_{c} = b_{c} |c| -1$.
Thus $\sGproj\L$ is $1$-Calabi-Yau by \cref{claim_34}. 
Since $\P_\L = \Eco_\L\cup\E_\L$ with $\Eco_\L\neq\E_\L$, it follows from \cref{claim_27,claim_20,claim_28,claim_41} that $\L$ is (non-G-semisimple) G-radical-square-zero, and 
\begin{align}
\sGproj\L=\add\left(\Omega_\L^{2}(\L/\rad\L) \oplus \Omega_\L^{3}(\L/\rad\L)\right)=\thick(\Omega_\L^{2}(\L/\rad\L))
\end{align}
and
\begin{align}
 \D_{\rm sg}(\L)=\add(\L/\rad\L\oplus\Sigma\hspace{1pt}\L/\rad\L)=\thick(\L/\rad\L).
\end{align}
Note that 
\begin{align}
\Omega_\L^{2}(\L/\rad\L) = a_2 a_1 a_2 a_1 a_2\L  
\quad\text{and}\quad
\Omega_\L^{3}(\L/\rad\L) = a_1 a_2\L,
\end{align}
where $a_2 a_1 a_2 a_1 a_2\L \cong a_1 a_2 a_1 a_2\L$ by \cite[Proposition 4.6]{Chen-Shen-Zhou_2018}. 
\end{example}




\appendix 
\section{Semisimple triangulated categories} 
\label{Semisimple_triangulated_categories} 

The aim of this appendix is to show that a Krull-Schmidt triangulated category $\T$ is semisimple if and only if $\mod\T$ is a semisimple abelian category.
We also show that these conditions are equivalent to the vanishing of the radical $\rad\T$ of $\T$. 
Note that Harada \cite[Theorem 1.2]{Harada_1970} and Street \cite[Proposition 11]{Street_1995}  characterized the vanishing of the radical of an artinian or noetherian abelian category and an additive category, respectively.

An abelian category is called {\it semisimple} if every short exact sequence splits. 
Every indecomposable object of a Krull-Schmidt semisimple abelian category is simple.
The following proposition characterizes Krull-Schmidt semisimple abelian categories. 
We learned the proof of (4) $\Rightarrow$ (2) from Osamu Iyama.
Before proceeding, we note that a sink morphism $u_X: S_X \to X$  (resp.~a source morphism $u^X: X \to S^X$) is zero if and only if $S_X=0$ (resp.~$S^X=0$); see \cite[Lemma 1.1]{ShipingLiu_2010}.

\begin{proposition}\label{claim_6}
The following conditions are equivalent for a Krull-Schmidt category $\C$.
\begin{enumerate}
\item $\C$ is a semisimple abelian category.
\item $\rad\C = 0$.
\item Every simple functor of $\Mod\C$ is a projective functor.
\item $\mod\C$ is a semisimple abelian category. 
\item For each $Z\in\ind\C$, the morphism $0\to Z$ is a sink morphism.
\item For each $X\in\ind\C$, the morphism $X\to 0$ is a source morphism.
\end{enumerate}
\end{proposition}
\begin{proof}
(1) $\Leftrightarrow$ (2)
Observe that $\rad\C=0$ if and only if any morphism in $\C$ is of the form $\left[\substack{\,1\,0\,\\[.5mm]\,0\,0\,}\right] : A \oplus B \to A \oplus C$ with $A,B,C\in\C$ (cf.~\cite[Proposition 3.1]{Rump_2022}).
The latter is equivalent to the condition that $\C$ admits a semisimple abelian structure.

(2) $\Leftrightarrow$ (3)
Follows from the fact that for each $X\in\ind\C$, the canonical surjection $\C(-, X)\to (\C/\rad\C)(-, X)$ is a projective cover of the simple functor $(\C/\rad\C)(-, X)$; see \Cref{Sink_and_source_morphisms}.

(2) $\Rightarrow$ (4) 
Assume (2) holds. 
Then $\mod\C$ is abelian, since $\C$ clearly has weak kernels.
Moreover, every $F \in \mod\C$ is projective.
Indeed, by \cite[Proposition A.1]{Krause_2012}, every $F\in\mod\C$ admits a minimal projective presentation $\C(-, X)\xrightarrow{(-, u)} \C(-, Y) \xrightarrow{\pi} F \to 0$ with $u \in\C(X, Y)$.
Since $\pi$ is a projective cover of $F$,  \cite[Proposition 3.10]{Krause_2015} implies that $(-, u)$ lies in 
\begin{align}
    \rad(\mod\C)(\C(-, X), \C(-, Y)) =\{(-, u):\C(-, X) \to \C(-, Y) \mid u \in \rad\C(X, Y)\}
    =0.
\end{align}
This shows that the projective cover $\pi: \C(-, Y) \to F$ is an isomorphism.
Therefore, $\mod\C$ is a semisimple abelian category.

(4) $\Rightarrow$ (2)
Take a non-zero morphism $u\in\C(X, Y)$ between indecomposable objects, and consider the induced exact sequence $\C(-, X)\xrightarrow{(-,u)}\C(-, Y) \to F\to 0$, where $F$ denotes the cokernel of $(-,u)$.
Since $F\in\mod\C$ is projective by assumption and $\Im(-,u)\neq 0$, we have $F=0$, equivalently, $u$ is an isomorphism. 
Therefore, $\rad\C(X, Y) = 0$.

We now prove (1) $\Leftrightarrow$ (5); the proof of (1) $\Leftrightarrow$ (6) is dual.

(1) $\Rightarrow$ (5)
Suppose $\rad\C=0$.
Then any non-zero morphism between indecomposable objects is an isomorphism.
This implies that $0\to Z$ is a sink morphism for each $Z\in \ind\C$.

(5) $\Rightarrow$ (1)
Assume (5) holds.
Then, for any $Z \in\ind\C$, the morphism $(-, 0) : \C(-, 0) \to \C(-, Z)$ induced by the sink morphism $0: 0 \to Z$ has a simple cokernel $(\C/\rad\C)(-, Z)$.
This forces $\rad\C(-, Z)$ to be zero.
\end{proof}

\begin{remark}\label{remark_1}
Combining \cref{claim_6}, \cite[Theorem 2.9 and Lemma 3.2]{Iyama-Nakaoka-Palu_2024} and \cite[Theorem 3.3]{Msapato_2025}, we see that extriangulated categories $\C$ satisfying $\rad\C=0$ are triangulated extriangulated categories.
\end{remark}

\begin{example}\label{example_7}
Let $D$ be a division ring.
Then $\mod D$ is a Krull-Schmidt semisimple abelian category such that $\ind\mod D =\{D\}$. 
Moreover, the bounded derived category $\D^{\rm b}(\mod D)$ of $\mod D$ is a Krull-Schmidt semisimple abelian category.
Indeed, we have $\ind\D^{\rm b}(\mod D) =\{\Sigma^i D\,|\,i \in \Z\}$ and 
\begin{align}
    \D^{\rm b}(\mod D)(\Sigma^iD, \Sigma^jD)
    \cong
    \begin{cases}
        D & \text{if } i = j,\\
        0 & \text{otherwise},
    \end{cases}
\end{align}
so that $\rad\D^{\rm b}(\mod D)=0$.
Further, for each $n \geq 1$, the orbit category $\D^{\rm b}(\mod D)/\Sigma^n$ is a Krull-Schmidt semisimple abelian category.
Indeed, we have $\ind\D^{\rm b}(\mod D)/\Sigma^n =\{\Sigma^iD\,|\,1 \leq i \leq n\}$ and 
\begin{align}
    (\D^{\rm b}(\mod D)/\Sigma^n)(\Sigma^iD, \Sigma^jD) = \bigoplus_{k\in\Z}\D^{\rm b}(\mod D)(\Sigma^iD, \Sigma^{j+kn}D)
    \cong 
    \begin{cases}
        D & \text{if } i = j,\\
        0 & \text{otherwise},
    \end{cases}
\end{align}
so that $\rad\D^{\rm b}(\mod D)/\Sigma^n =0$.
\end{example}

For a semisimple abelian category $\mathcal{A}$ and its
auto-equivalence $\Sigma$, there is a unique triangulated structure on $\mathcal{A}$ with $\Sigma$ the shift functor; see \cite[Lemma 3.4]{X-WChen_2011}.
The resulting triangulated category is denoted by $(\mathcal{A}, \Sigma)$.

\begin{definition}
A triangulated category is called {\it semisimple} if it is triangle equivalent to $(\mathcal{A}, \Sigma)$ for some semisimple abelian category $\mathcal{A}$ and some auto-equivalence $\Sigma:\mathcal{A}\to \A$.
\end{definition}

Recall that a quiver is {\it discrete} if it has no arrow.
The following theorem characterizes Krull-Schmidt semisimple triangulated categories.
Note that it is implicitly contained in \cite[Proposition 6.2]{ShipingLiu_2010}.

\begin{theorem}\label{claim_5}
The following conditions are equivalent for a Krull-Schmidt triangulated category $\T$.
\begin{enumerate}
\item $\T$ is a semisimple triangulated category.
\item $\rad\T=0$.
\item For each $Z \in \ind\T$, the  triangle $\Sigma^{-1}Z \to 0 \to Z \xrightarrow{1} Z$ is Auslander-Reiten.
\item $\T$ has Auslander-Reiten triangles, and the Auslander-Reiten quiver of $\T$ is discrete.   
\end{enumerate}
\end{theorem}
\begin{proof}

(1) $\Leftrightarrow$ (2) $\Leftrightarrow$ (3)
Follows from \cref{claim_6}.
Note that (2) $\Rightarrow$ (1) is obtained from the observation that  the identity functor $1_\T: \T \to \T$ induces a triangle equivalence between $\T$ and the semisimple triangulated category $(\T, \Sigma)$.

(2)+(3) $\Rightarrow$ (4)
Assume (2) and (3) hold.
Then $\T$ has Auslander-Reiten triangles by (3).
Moreover, since there is no irreducible morphism in $\T$ by (2), the Auslander-Reiten quiver of $\T$ is discrete.

(4) $\Rightarrow$ (3)
Assume (4) holds.
Then $d_{X,Y} = 0$ for all $X, Y \in\ind\T$, where $d_{X,Y}$ is defined in \eqref{eq_2}.
This implies that all Auslander-Reiten triangles in $\T$ have zero middle term.
\end{proof}

\begin{example}\label{example_5}
Let $D$ be a division ring and $n\in\Z$.
By \cref{example_7} and \cref{claim_5}, the semisimple abelian orbit category $\D^{\rm b}(\mod D)/\Sigma^n$ and the auto-equivalence $\underline{\Sigma}$ of $\D^{\rm b}(\mod D)/\Sigma^n$ induced by the shift functor $\Sigma$ of $\D^{\rm b}(\mod D)$ form a semisimple triangulated category. 
Assume $n\geq 1$, and let $D^n$ be the direct product of $n$ copies of $D$.
Then $\mod D^n$ is a Krull-Schmidt semisimple abelian category with $|\ind\mod D^n| = n$.
Define an algebra automorphism $\sigma :D^n \to D^n$ by $\sigma(x_1, \ldots, x_n) = (x_2, \ldots, x_n, x_1)$ and denote by $\sigma^*$ the automorphism of $\mod D^n$ twisting the right $D^n$-action by $\sigma :D^n \to D^n$.
Then the two semisimple triangulated categories $(\mod D^n, \sigma^*)$ and $(\D^{\rm b}(\mod D)/\Sigma^n, \underline{\Sigma})$ are triangle equivalent.
Indeed, let $C:=\bigoplus_{1\leq i\leq n}\Sigma^i D$.
Since 
\begin{align}
    (\D^{\rm b}(\mod D)/\Sigma^n)(C, C) \cong \prod_{1\leq i\leq n}(\D^{\rm b}(\mod D)/\Sigma^n)(\Sigma^iD, \Sigma^iD)\cong D^n
\end{align}
as rings, there are equivalences 
\begin{align}
    \D^{\rm b}(\mod D)/\Sigma^n 
    \xrightarrow{\sim} \proj(\D^{\rm b}(\mod D)/\Sigma^n)(C, C)
    \xrightarrow{\sim} \proj D^n =\mod D^n
\end{align}
of Krull-Schmidt semisimple abelian categories, where the first equivalence follows from \cite[Proposition 2.7]{Auslander_1974_I}.
Since the obtained equivalence $F:\D^{\rm b}(\mod D)/\Sigma^n\to \mod D^n$ is given by $(\D^{\rm b}(\mod D)/\Sigma^n)(C, -)$, it follows that $F \underline{\Sigma} \cong \sigma^* F$, which yields the desired triangle equivalence  
\begin{align}
    (\D^{\rm b}(\mod D)/\Sigma^n, \underline{\Sigma}) \xrightarrow{\sim} (\mod D^n, \sigma^*).
\end{align}
\end{example}

There is another example of a semisimple triangulated category.
It was proved by Chen  \cite[Theorem 3.8]{X-WChen_2011} that the singularity category $\D_{\rm sg}(\L)$ of an artin algebra  $\L$ with radical square zero is a semisimple triangulated category, 
where the corresponding semisimple abelian category and auto-equivalence are both explicitly described.

We now specialize \cref{claim_40} to semisimple triangulated categories.
Note that the definition of path-closed subsets is given preceding the proposition. 
We also note that the first statement of the following proposition is implicitly contained in \cite[Lemma 4.3]{Chen-Ringel_2018}.

\begin{proposition}\label{claim_37}
Let $\T$ be a Krull-Schmidt semisimple triangulated category.
\begin{enumerate}
\item For each path-closed subset $\U$ of $\ind\T$, its additive closure $\T_\U=\add\U$ in $\T$ is a connected semisimple thick subcategory of $\T$ such that for any $X\in\U$, there is a triangle equivalence
\begin{align}
    \T_\U = \thick_\T(X) \cong (\D^{\rm b}(\mod \T(X, X))/\Sigma^{n}, \underline{\Sigma})
\end{align} 
for some unique $n\in\Z$, where $\underline{\Sigma}$ is the auto-equivalence of $\D^{\rm b}(\mod \T(X, X))/\Sigma^{n}$ induced by the shift functor $\Sigma$ of $\D^{\rm b}(\mod \T(X, X))$.
\item For any two distinct path-closed subsets $\U, \U^\prime$ of $\ind\T$, we have 
\begin{align}
    \T(\T_{\U}, \T_{\U^\prime}) =0= \T(\T_{\U^\prime}, \T_{\U}).
\end{align}
In particular, if $\ind\T$ has a finite number of path-closed subsets, then there is a decomposition $\T\,=\,\prod_{\U}\T_\U$ into a finite product of connected thick subcategories, where $\U$ runs over the path-closed subsets of $\ind\T$.
\end{enumerate}
\end{proposition}
\begin{proof}
(1)
Let $\U$ be a path-closed subset of $\ind\T$ and $X\in\U$.
Since $S_X=0$ by \cref{claim_5}, it follows from \cref{claim_40} that $\T_\U=\thick_\T(S_X\oplus X)=\thick_\T(X)$.
Since there is no non-zero non-isomorphism in $\T_\G$ by assumption, we conclude from \cite[Lemma 4.3]{Chen-Ringel_2018} and \cref{example_5} that $\T_\U$ is triangle equivalent to $(\D^{\rm b}(\mod \T(X, X))/\Sigma^{n}, \underline{\Sigma})$ for some unique $n\in\Z$, as desired.

(2) Follows from \cref{claim_40}.
\end{proof}

We end this section with the following observation.

\begin{proposition}\label{claim_12}
Let $\T$ be a $k$-linear Hom-finite Krull-Schmidt semisimple triangulated category.
Assume further that $\T(X, X) \cong k$ for all $X\in\ind\T$.
Then $\T$ is $0$-Calabi-Yau.
\end{proposition}
\begin{proof} 
By \cref{claim_5}, the action of the Serre functor $\nu$ for $\T$ on the indecomposable objects is the identity.
We show that $\nu$ acts on the morphisms between indecomposable objects by the identity.
For $X, Y \in\ind\T$ with $\T(X, Y)\neq 0$, we fix a $k$-basis element $f_{X, Y}$ of the one dimensional vector space $\T(X, Y)$ in a way that $f_{Y, X} = f_{X, Y}^{-1}$ and  $f_{X, X} = \id_X$. 
Consider the pairing $\langle-,-\rangle:\T(X, Y)\times\T(Y, X) \to \T(X, X) \xrightarrow[\sim]{\,\eta_X\,} k$ given by $\langle f, g \rangle = \eta_X(gf)$ for  $(f, g) \in \T(X, Y)\times\T(Y, X)$, where $\eta_X$ is the $k$-linear isomorphism sending $f_{X, X}$ to $1$.
Then $\langle-,-\rangle$ is non-degenerate (see \cite[Lemma 6.3.1]{Krause_book}), and the isomorphisms
\begin{align}
    \eta_{X, Y}:\T(X, Y) \xrightarrow{\,\sim\,} D\T(Y, X) 
    \quad\text{and}\quad
    \eta_{X, Y}^\prime:\T(Y, X) \xrightarrow{\,\sim\,} D\T(X, Y) 
\end{align}
induced by the pairing $\langle-,-\rangle$ send $f_{X, Y}$ and $f_{Y, X}$ to $\langle f_{X, Y},-\rangle$ and $\langle -, f_{Y, X}\rangle$, respectively.
A direct calculation shows that the composition 
\begin{align}
    \T(X, Y) \xrightarrow[\sim]{\,\varepsilon\,} DD\T(X, Y) \xrightarrow[\sim]{\,D\eta_{X, Y}^\prime\,} D\T(Y, X)  
\end{align}
sends $f_{X, Y}$ to $\langle f_{X, Y},-\rangle$, where $\varepsilon$ is the canonical isomorphism given by $\varepsilon(f_{X, Y})(\varphi) = \varphi(f_{X, Y})$ for $\varphi \in D\T(X, Y)$.
Therefore, the composition  
\begin{align}
    \T(X, Y) \xrightarrow[\sim]{\,\eta_{X, Y}\,} D\T(Y, X) \xrightarrow[\sim]{\left(D\eta_{X, Y}^\prime\,\circ\,\varepsilon\right)^{-1}} \T(X, Y)
\end{align}
is the identity.
This composition is  the isomorphism induced by the Serre functor $\nu $ (see \cite[Lemma 6.3.1]{Krause_book} again).
\end{proof}




\section*{Acknowledgments} 
The second author would like to thank Rasool Hafezi and Osamu Iyama for valuable discussions.

\bibliographystyle{plain}
\bibliography{ref}

\end{document}